\documentclass[11pt]{article}
\usepackage{amsfonts, amsmath, amssymb}
\usepackage{graphicx, amsthm, color}
\newtheorem{theorem}{Theorem}[section]
\newtheorem{lemma}[theorem]{Lemma}

\newtheorem{corollary}[theorem]{Corollary}

\newtheorem{definition}[theorem]{Definition}
\theoremstyle{definition}

\theoremstyle{remark}
\newtheorem{remark}[theorem]{Remark}
\numberwithin{equation}{section}

\definecolor{ForestGreenWeb}{rgb}{0.13, 0.55, 0.13}
\begin{document}

\begin{center}
{\LARGE All convex bodies are in the subdifferential of some\bigskip \\
everywhere differentiable locally Lipschitz function}

\vspace{1cm}

{\Large \textsc{Aris Daniilidis, Robert Deville, Sebasti{\'{a}}n
Tapia-Garc{\'{\i}}a}}
\end{center}

\bigskip

\noindent\textbf{Abstract.} We construct a differentiable locally Lipschitz function $f$ in
$\mathbb{R}^{N}$ with the property that for every convex body $K\subset \mathbb{R}^N$ there exists $\bar x \in \mathbb{R}^N$ 
such that $K$ coincides with the set $\partial_L f(\bar x)$ of limits of derivatives $\{Df(x_n)\}_{n\geq 1}$ of sequences $\{x_n\}_{n\geq 1}$ converging to~$\bar x$.
The technique can be further refined to recover all compact connected subsets with nonempty interior, disclosing an important difference between differentiable and continuously differentiable functions. It stems out from our approach that the class of these pathological functions contains an infinite dimensional vector space and is dense in the space of all locally Lipschitz functions for the uniform convergence.

\bigskip


\vspace{0.6cm}

\textit{Secondary} 26A24, 26A16, 26A27

\medskip

\noindent\rule{17cm}{0.8mm}\\ 
\vspace{-0.7cm}

\tableofcontents 
\vspace{2mm}
\noindent\rule{17cm}{0.8 mm} 


\section{Introduction}

Given a nonempty open subset $\mathcal{U}$ of a Euclidean space $\mathbb{R}^{N}$, a function $f:\mathcal{U}\rightarrow\mathbb{R}$ is called Lipschitz if
there exists a constant $L>0$ such that
\begin{equation}
|f(x)-f(y)|\,\leq\,L\,\Vert x-y\Vert,\quad\text{for all }x,y\in\mathcal{U}.
\label{Lip}
\end{equation}
We denote by $\|f\|_{\mathrm{Lip}}$ the infimum of the above constants, so
that:
\begin{equation}
\|f\|_{\mathrm{Lip}}=\sup_{x,y\in\mathcal{U}\,,\,x\neq y}\,\frac
{|f(x)-f(y)|}{\Vert x-y\Vert}. \label{LipConst}
\end{equation}
In what follows, we call a function $k$-Lipschitz if $\|f\|_{\mathrm{Lip}}\leq
k$, where $k\geq0.$ We also call a function locally Lipschitz, if around any
point $x_{0}$ of its domain, there exists $k>0$ and a neighborhood $\mathcal{V}$ of
$x_{0}$ such that the function $f$ is $k$-Lipschitz on $\mathcal{V}$.
\smallskip

According to the Rademacher theorem, every locally Lipschitz function is
differentiable almost everywhere (see \cite[Chapter~9]{BLbook}\textit{ e.g.}).
If $\mathcal{N}$ is any null subset of $\mathcal{U}\subset\mathbb{R}^N$, then denoting by
$\mathcal{D}_{f}$ the set of points of differentiability of $f$ and by $Df(x)$
the derivative of $f$ at a point $x\in\mathcal{D}_{f}$, the Clarke subdifferential at
$x\in\mathcal{U}$ is given by the following formula (see \cite[Chapter~2]{Clarke}):
\begin{equation}
\partial f(x)=\mathrm{conv}\left\{  \lim_{x_{n}\rightarrow x}\ Df(x_{n}
)\,:\,\{x_{n}\}\subseteq\mathcal{D}_{f}\diagdown\mathcal{N} \right\}  ,
\label{clarke}
\end{equation}
where $\mathrm{conv}(A)$ stands for the convex envelope of a set $A$.
It follows that the above definition is independent of the choice of
$\mathcal{N}$ and that $\partial f(x)$ is a nonempty convex compact subset of
the closed dual ball $\overline{B}(0,\|f\|_{\mathrm{Lip}})$ containing the
derivative $Df(x)$, whenever this latter exists. \smallskip

{The Clarke subdifferential admits an alternative description based on \textit{Fréchet} subgradients, without explicit use of derivatives or the above null set. We recall that $x^{\ast}\in\mathbb{R}^N$ is a Fréchet subgradient of $f$ at $x$ (and denote $x^{\ast}\in \widehat{\partial}f(x)$) if $x^{\ast}=\nabla \phi(x)$ for some $\mathcal{C}^1$-smooth function $\phi\leq f$ with $\phi(x)=f(x)$. Then we say that $p\in \mathbb{R}^{N}$ is a \textit{limiting} subgradient of $f$ at $x$, and denote $p \in\partial_L f(x)$, if there exists a sequence $\{ (x_n,x^{\ast}_n) \}_n$ in $\mathbb{R}^N \times \mathbb{R}^N$ with $x_n^{\ast}\in \widehat{\partial}f(x_n)$ such that $\mathrm{lim}_{n \to \infty} x_n=x$ and $\mathrm{lim}_{n \to \infty} x_n^{\ast}=p$. The Clarke subdifferential can then be defined as the convex envelope of the limiting subdifferential, that is, $\partial f(x)=\mathrm{conv}\,\{\partial_L f(x)\}$, for every locally Lipschitz function $f$. Therefore, $\partial_L f(x) \subset \partial f(x)$. Notice that if $f$ is everywhere differentiable, the limiting subdifferential is given by the formula:
\begin{equation}\label{eq-boris}
\partial_L f(x) :=\left\{  \lim_{x_{n}\rightarrow x}\ Df(x_{n}) \right\}  
\end{equation} 
and if} $f$ is $\mathcal{C}^{1} $-smooth, we have $\partial f(x)=\partial_L f(x) = \{Df(x)\}$, for all
$x\in\mathcal{U}$. In fact, {for a Lipschitz function $f$,} $\partial f(x)$ reduces to a singleton if and only if $f$ is strictly differentiable at $x$ {\cite[Proposition 2.2.4]{Clarke}}.

\smallskip

Notice that $\|\cdot\|_{\mathrm{Lip}}$ is a seminorm in the vector space
$\mathrm{Lip}(\mathcal{U})$ of all real-valued Lipschitz functions on
$\mathcal{U}$ and becomes a norm in the subspace $\mathrm{Lip}_{x_{0}
}(\mathcal{U})$ of those functions that vanish at some (arbitrarily chosen)
prescribed point $x_{0} \in\mathcal{U}$. In particular, $(\mathrm{Lip}_{x_{0}}
(\mathcal{U}), \|\cdot\|\mathrm{_{\mathrm{Lip}}})$ is a Banach space (known
also as the dual space of the free space of $\mathrm{\mathcal{U}}$).
Alternatively, setting $\|\cdot\|_{\mathrm{L}}:=\|\cdot\|_{\infty}
+\|\cdot\|_{\mathrm{Lip}}$ and denoting by $\mathcal{L}^\infty(\mathcal{U})$ the set of bounded
functions on $\mathcal{U}$, the normed space $(\mathrm{Lip}(\mathcal{U}
)\cap\mathcal{L}^\infty(\mathcal{U}),\|\cdot\|_{\mathrm{L}})$ is also complete.
 \smallskip

If the set $\mathcal{U}$ is bounded, one can also consider the norm
$\|f\|_{\infty}:=\sup\,\left\{  |f(x)|:x\in\mathcal{U}\right\}  $ of uniform
convergence. In this case $\left(  \mathrm{Lip}(\mathcal{U}),\|\cdot
\|_{\infty}\right)  $ is not complete (in fact, it is dense in the Banach space
$\left(  \mathcal{C}_{b}(\mathcal{U}),\|\cdot\|_{\infty}\right)  $ of bounded
continuous functions). However, one can remedy this lack of completeness by
considering the set $\mathrm{Lip}^{[k]}(\mathcal{U})$ of Lipschitz
continuous functions with Lipschitz constant $\|f\|_{\mathrm{Lip}}\leq k$.
This set is a complete metric space under the distance of uniform
convergence $d_{\infty}(f,g):=\|f-g\|_{\infty}$. In this setting (where the vector structure is of course lost) {and assuming that $\mathcal{U}$ is convex}, a standard application of Baire's category theorem has been used by J. Borwein and X.~Wang (see \cite{BW2000, BW2003} \textit{e.g.}) to {establish} that the set of Lipschitz functions with maximal Clarke subdifferential (that is, $\partial
f(x)\equiv\overline{B}(0,k)$ for all $x\in\mathcal{U}$) is residual in
$\mathrm{Lip}^{[k]}(\mathcal{U})$. Therefore, a generic Lipschitz function in $\mathrm{Lip}^{[k]}(\mathcal{U)}$ has Lipschitz constant equal to $k$ and \textit{saturates} its Clarke subdifferential at every point. The first explicit construction of a Clarke saturated function was given in \cite{Lebourg} (in dimension one) and in \cite{BS2010} (in higher dimensions).

\smallskip

The aforementioned result of J. Borwein and X. Wang {underlines} the fact that uniform convergence does not entail any control on derivatives and local oscillations. The genericity is thus tightly related to the $d_{\infty}
$-topology: one easily sees that the set of Clarke-saturated functions (that
is, functions whose subdifferential is identically equal to the closed ball
$\overline{B}(0,\|f\|_{\mathrm{Lip}})$) cannot be dense for the (more adequate)
distance $d_{\mathrm{Lip}}(f,g)=\|f-g\|_{\mathrm{Lip}}$ given by the Lipschitz
norm. Still, in \cite{DF2019} it was established that the set of Clarke saturated
functions is \textit{spaceable} in $(\mathrm{Lip}(\mathcal{U}),\|\cdot
\|_{\mathrm{L}}),$ that is, it contains a closed infinite dimensional
subspace (see \cite{AGS2005, EGS2014} for a discussion about
spaceability). The construction of this infinite dimensional subspace of
Clarke saturated functions is explicit, but the result requires working in
$\ell_{1}^{N}$ (rather than in the usual Euclidean space $\mathbb{R}^{N}$).

\smallskip

Let us mention for completeness that important subclasses of Lipschitz
functions, such as semialgebraic (more generally, Whitney stratifiable) or
finite selections of $\mathcal{C}^{N}$-smooth functions have small Clarke
subdifferentials: they often reduce to a singleton and the (generalized)
critical values satisfy the conclusion of the Morse-Sard theorem, see
\cite[Corollary~5(ii)]{BDLS} and \cite[Theorem~5]{BDD} respectively. On the
other hand, every point of a Clarke-saturated Lipschitz function is {(Clarke) critical, since 
$0\in \partial f(x)\equiv \overline{B}(0,\|f\|_L)$.}
Other pathological situations have also been detected in~\cite{DD2020} where
the authors constructed examples of Lipschitz continuous functions with finite
Clarke critical values, but with pathological subgradient dynamics both in
continuous and discrete time: the iterates generate bounded trajectories that fail to detect any Clarke critical point of the function. Finally, in~\cite{BMW1997} the authors constructed locally Lipschitz functions whose subdifferential assumes a prescribed set of values. \smallskip

In this work we establish the following result for the range of the Clarke subdifferential. (The term \textit{convex body} employed below will refer to a compact convex set with nonempty interior.)
\begin{itemize}
\item There exists a compactly supported, differentiable $1$-Lipschitz function $f:\mathbb{R}^{N}\rightarrow \mathbb{R}$ whose Clarke
subdifferential contains all convex bodies of the closed unit ball.
\end{itemize}
\noindent The construction is different for $N=1$ (Theorem~\ref{1D}) and for $N\geq 2$ (Theorem~\ref{theo: function RN}). In the first case, the function $f$ is also \textit{subdifferentially exhaustive} (see Definition~\ref{def_exh}), that is, its Clarke subdifferential takes all of its possible values. In both cases, $N=1$ and $N\geq 2$, the construction reveals that the set of all such functions whose support is contained in an open bounded set $\mathcal{U}$ of $\mathbb{R}^{N}$ is \textit{spaceable} in $(\mathrm{Lip}(\mathcal{U}),\Vert \cdot \Vert _{\mathrm{L}})$ and dense in $(\mathrm{Lip}^{[1]}(\mathcal{U}),d_{\infty })$, see Remark~\ref{dense}~(v),(vi) and Subsection~\ref{ss-3.2}. 

\smallskip

By enhancing the techniques employed in Subsection~\ref{ss-3.2} we obtain, in Subsection~\ref{ss-tapia}, a more general result (Theorem~\ref{jortega}), that recovers all compact connected subsets of $\mathbb{R}^N$ with nonempty interior (not only the convex bodies). The construction requires $N\geq 2$ (but for $N=1$ the two notions coincide anyway). The general result reads as follows:
\begin{itemize}
\item There exists a compactly supported, differentiable function $f:\mathbb{R}^{N}\rightarrow \mathbb{R}$ 
whose limiting subdifferential contains all compact connected subsets of $\mathbb{R}^{N}$ with nonempty interior.
\end{itemize}

\section{Prerequisites.}

We recall that the term Polish space refers to any separable topological space, whose topology can be
metrizable in a way that the resulting metric space is complete. We denote by $\Delta:=\{0,1\}^{\mathbb{N}}$ the Cantor set and recall
that every uncountable Polish space contains a homeomorphic copy of~$\Delta$, see \cite[Corollary~6.5]{Kechris}. \smallskip\newline 
In this work, we consider the Euclidean space $\mathbb{R}^N$, $N\geq 1$ and
denote by $B(0,r)$ (respectively, $\overline{B}(0,r)$) the open (respectively, closed) ball 
centered at $x\in\mathbb{R}^N$ with radius $r>0$. \smallskip\newline
Given a nonempty convex compact subset $C$ of $\mathbb{R}^{N}$ we set: 
\begin{equation}
\mathcal{F}_{C}:=\{K\subset C:\,K\neq \emptyset ,\text{ compact}\}.
\label{eq:F}
\end{equation}
It is known that $\mathcal{F}_{C}$ is a compact metric space for the
Hausdorff distance 
\begin{equation}
D_{\mathrm{H}}(K_{1},K_{2}):=\max \,\left\{ \underset{x\in K_{1}}{\sup }d(x,K_{2}),\,\underset{x\in K_{2}}{\sup }d(x,K_{1})\right\}   \label{eq:DH}
\end{equation}
where $d(x,A):=\inf \left\{ \Vert x-a\Vert :a\in A\right\} $ for every $A\subset \mathbb{R}^{N}$. We further set
\begin{equation}
\mathcal{K}_{C}:=\{K\subset C:\,K\neq \emptyset ,\text{ compact convex}\}.
\label{eq:G}
\end{equation}
Notice that $\mathcal{K}_{C}$ is a closed subset of $\mathcal{F}_{C}$ under
the Hausdorff distance, therefore $(\mathcal{K}_{C},D_{\mathrm{H}})$ is also
a compact metric space.
\smallskip

\noindent In what follows we denote by $\mathcal{L}_{N}$ the Lebesgue measure
on $\mathbb{R}^{N}$. Given an integrable function $f:\mathbb{R}^{N}
\rightarrow\mathbb{R}$, we say that a point $x$ is a \emph{Lebesgue point} of
$f$ if
\[
\lim_{r\searrow0^{+}}\frac{1}{\mathcal{L}_{N}(B(x,r))}{\displaystyle\int\limits_{B(x,r)}}
|f(y)-f(x)|\,dy=0
\]
Therefore, a Lebesgue point is a point where $f$ does not oscillate in an
average sense, see~\cite[§1.7]{EG-book}. It is known that the set of Lebesgue points
of every integrable function $f$ is of full-measure. In particular, for
\emph{a.e.} $x\in\mathbb{R}^{N}$ it holds
\[
\left\vert \, f(x) - \left( \frac{1}{\mathcal{L}_{N}(B(x,r))} {\displaystyle\int\limits_{B(x,r)}}
f(y)dy\right)\,\right\vert \underset{r\rightarrow0}{\longrightarrow}
0\qquad\text{(Lebesgue differentiation theorem)}
\]
Let us further recall the \textit{interval splitting property} for subsets of the
real line.
\begin{definition}
[splitting property]\label{def-everywhere}$\mathrm{(i).}$ A set $A\subset\mathbb{R}$ is
called everywhere positive-measured, if it intersects any nontrivial interval
in a set of positive Lebesgue measure.\smallskip\newline$\mathrm{(ii).}$ We say that
$A$ has the splitting property for the family of intervals of $\mathbb{R}$ if
both $A$ and $\mathbb{R}\diagdown A$ are everywhere positive-measured.
\end{definition}

The following lemma goes back to Bruckner \cite{Bruckner-1967} (see also
\cite[Lemma~4.1]{Wang-2004}).

\begin{lemma}
[countable splitting partition]\label{CounPart}There exists a countable
partition $\{A_{k}\}_{k\in\mathbb{N}}$ of $\mathbb{R}$, each of which splits
the family of intervals.
\end{lemma}

\smallskip

Let us now recall that given a nonempty open subset $\mathcal{U}$ of
$\mathbb{R}^{N}$ and a $k$-Lipschitz function $f:\mathcal{U}\rightarrow
\mathbb{R}$, the Clarke subdifferential operator $\partial f:\mathcal{U}
\rightrightarrows\mathbb{R}^{N}$ has closed graph and nonempty convex compact
values (in particular, $\partial f(x)\subset\overline{B}(0,k)$ for every
$x\in\mathcal{U}$). We also recall that $\partial f$ is an \emph{upper
semicontinuous} multivalued operator, in the sense that for every
$\varepsilon>0$ and $x\in\mathcal{U}$ there exists $\delta>0$ such that for
all $y\in B(x,\delta)\cap\mathcal{U}$ it holds $\partial f(y)\subset\partial
f(x)+B(0,\varepsilon).$\smallskip\newline In what follows, $\mathcal{U}$ will denote a nonempty
open subset of $\mathbb{R}^{N}$. We recall from \cite{DF2019} the
following definition.

\begin{definition}
[subdifferential saturation]\label{def_sat}  A Lipschitz function $f:\mathcal{U}
\rightarrow\mathbb{R}$ is called \textit{Clarke saturated} if for every
$x\in\mathcal{U}$ we have $\partial f(x)=\overline{B}(0,\|f\|_{\mathrm{Lip}})$. 
\end{definition}
\noindent Therefore, a Lipschitz function $f$ with $\|f\|_{\mathrm{Lip}}=1$ is Clarke saturated if and only if its Clarke
subdifferential at any point is equal to the unit ball of $\mathbb{R}^{N}$. \smallskip

\noindent We shall further use the following terminology.

\begin{definition}
[subdifferential exhaustiveness]\label{def_exh} A Lipschitz function $f:\mathcal{U}
\rightarrow\mathbb{R}$ is called Clarke exhaustive (respectively, almost
exhaustive) if for any nonempty closed convex subset $K$ (respectively, of
nonempty interior) of the ball $\overline{B}(0,\|f\|_{\mathrm{Lip}})$, there exists
$x\in\mathcal{U}$ such that $\partial f(x)=K$.
\end{definition}


\section{Main results}

In this section we are going to construct an everywhere differentiable function in $\mathbb{R}^N$ with
bounded derivatives (thus, in particular, a Lipschitz continuous function) whose Clarke
subdifferential is almost exhaustive. This yield the result announced in the title of the paper.\smallskip\newline 
The construction requires at least two dimensions (that is, $N\geq 2$), but the result is also true for $N=1$ through 
a different construction which will be treated first. Moreover, in the $1$-dimensional case the constructed function 
turns out to be Clarke-exhaustive, that is, the subdifferential is surjective (assuming all of its possible values). \smallskip\newline
Since we deal with functions which are everywhere differentiable, the result is rather unexpected, taking into account that the derivative 
is a Baire--1 function (therefore, gene\-rically continuous) and the Clarke subdifferential of a strictly differentiable function 
(thus, \textit{a fortiori}, of a $C^{1}$-function) is singleton everywhere.\smallskip\newline
As a matter of fact, our results also hold for the (smaller) limiting subdifferential, see forthcoming Remark~\ref{dense}(i) (for $N=1$) and Remark~\ref{devil} (for $N\geq 2$).
A further refinement will be performed in Section~\ref{ss-tapia} where we eventually show that there exists a differentiable, locally Lipschitz function such that
every compact connected subset of $\mathbb{R}^{N}$ with nonempty interior appears in the range of its limiting subdifferential.

\subsection{Subdifferentially exhaustive differentiable functions in $\mathbb{R}$}

Let $f:(0,1)\rightarrow\mathbb{R}$ be $1$-Lipschitz. Then for every
$x\in(0,1)$, the subdifferential $\partial f(x)$ is a nonempty closed
subinterval of $[-1,1]$ (possibly reducing to a singleton). We shall need the
following notation:
\begin{equation}
\mathbb{T}^{+}=\big\{(a,b)\in\mathbb{R}^{2}:\quad0\leq a\leq b\leq
1\big\} \label{q:T+}
\end{equation}
Let us start with the following essentially known result.

\begin{lemma}
\label{Cantor0} There exists a continuous surjective curve $\gamma_{1}:[0,1]\rightarrow\mathbb{T}^{+}$ such that ${\gamma_1((0,1))=\mathbb{T}^{+}}$.
\end{lemma}

\noindent\textit{Proof.} It is well-known that there exists a continuous
surjective curve $\gamma_{0}:[0,1]\rightarrow\lbrack0,1]\times\lbrack0,1]$.
This map is called a Peano curve, see \cite{P1890}. The function
$\varphi:\mathbb{R}^{2}\rightarrow\mathbb{R}^{2}$ defined by
\[
\varphi(a,b)=\bigl(\min\{a,b\},\max\{a,b\}\bigr)
\]
is continuous and maps $[0,1]\times\lbrack0,1]$ onto $\mathbb{T}^{+}$. Thus
the function $\gamma_{1}=\varphi\circ\gamma_{0}$ satisfies the assertion of
the statement. \hfill$\Box$

\bigskip

\subsubsection{An easy nonsmooth example}

Let us first provide a straightforward construction of a $1$-Lipschitz
Clarke exhaustive function (omitting momentarily the additional requirement of
being everywhere differentiable).

\begin{theorem}
[exhaustive Lipschitz function in $\mathbb{R}$]\label{nonsmooth} There exists
a Lipschitz function ${f:[0,1]\rightarrow\mathbb{R}}$ with $\|f\|_{\mathrm{Lip}}=1$ such that 
for every nonempty closed interval $[a,b]\subset\lbrack-1,1]$, there exists $x\in(0,1)$
such that $\partial f(x)=[a,b]$, that is, $\partial f([0,1])=\mathcal{K}_{[-1,1]}$.
\end{theorem}

\noindent\textit{Proof.}  Let $\gamma_{1}(t)=\bigl(a(t),b(t)\bigr)$, with $t\in\lbrack0,1]$,
be the continuous curve given by Lemma~\ref{Cantor0} and let $A\subset (0,1)$ be a measurable set
which splits the family of nonempty open intervals of $[0,1]$ (\textit{c.f.} Definition~\ref{def-everywhere}(ii)). 
The required function $f$ is explicitly defined as follows:
\[
f(x)=\int_{0}^{x}\big[a(t)\mathbf{1}_{A}(t)+b(t)\mathbf{1}_{[0,1]\backslash
A}(t)\big]dt.
\]
Indeed, let us prove that for every $x\in [0,1]$, we have $\partial
f(x)=[a(x),b(x)]$. \smallskip\newline To this end, let us first consider a
Lebesgue point $t\in A$ of the function $\mathbf{1}_{A}$. Since $a$ is
continuous, we have that $f^{\prime}(t)$ exists and $f^{\prime}(t)=a(t)$.
Similarly, if $s\in(0,1)\setminus A$ is a Lebesgue point of the function $\mathbf{1}_{[0,1]\backslash A}$, then $f^{\prime}(s)$ exists and $f^{\prime}(s)=b(s)$.
Fix now $x\in (0,1)$ (arbitrarily chosen). Since any open interval
containing $x$ meets the sets $A$ and $[0,1]\backslash A$ on a set of positive
measure, we deduce that $a(x)\in\partial f(x)$ and $b(x)\in\partial f(x),$
yielding $[a(x),b(x)]\subset\partial f(x)$.\smallskip\newline 
To establish the other inclusion, let us fix $\varepsilon>0$ and 
\[
\mathcal{N}:=\{x \in [0,1]:\, x \text{ is not a Lebesgue point for } f \}.
\]
Since the functions $a,b$ are
continuous, there exists $\delta>0$ such that $|a(t)-a(x)|\leq\varepsilon$ and
$|b(t)-b(x)|\leq\varepsilon$, for all $t\in(x-\delta,x+\delta)$. It follows
that $\partial f(x)\subset\lbrack a(x)-\varepsilon,b(x)+\varepsilon]$. Since
$\varepsilon>0$ is arbitrarily chosen, we deduce $\partial f(x)\subset\lbrack
a(x),b(x)]$ and consequently, equality holds.\medskip

Recalling that $\gamma_{1}$ satisfies Lemma~\ref{Cantor0}, for every nonempty
closed interval $[a,b]\subset\lbrack0,1]$, there exists $x\in(0,1)$ such that
$\gamma_{1}(x)=(a,b)\in\mathbb{T}^{+}\subset\mathbb{R}^{2}$, and
consequently, $\partial f(x)=[a,b]$. Replacing~$f$ by the function
\[
\widetilde{f}(x):=2f(x)-x,\qquad\text{for all }x\in\lbrack0,1],
\]
we obtain a function $\widetilde{f}$ which is also $1$-Lipschitz: Indeed, notice that $\widetilde{f}^{\prime}(x)=2f^{\prime}(x)-1\in\lbrack-1,1]$ whenever $f^{\prime}(x)$ exists. It follows directly that
\[
\partial\widetilde{f}((0,1))=\big\{[a,b];\,-1\leq a\leq b\leq
1\big\}=\mathcal{K}_{[-1,1]}.
\]
The proof is complete.\hfill$\Box$

\bigskip

\begin{remark}
Notice that the set of bounded $1$-Lipschitz Clarke exhaustive functions in $\mathbb{R}$
cannot be $d_{\infty}$-residual in the (complete) metric space $\left(
\mathrm{Lip}^{[1]}(\mathbb{R}),d_{\infty}\right)$ of all bounded
$1$-Lipschitz functions in $\mathbb{R}$, since it shares with the set
of Clarke-saturated functions (which is known to be $d_{\infty}$-residual, see~\cite{BW2000}) only the null function $f\equiv 0$.
However, we shall see later (Remark~\ref{dense}(iv)) that the set of bounded $1$-Lipschitz functions in $\mathbb{R}$ which are Clarke exhaustive is dense in $\left(\mathrm{Lip}^{[1]}(\mathbb{R}),d_{\infty}\right)$.
\end{remark}


\subsubsection{An involved construction ensuring differentiability}

We shall now enhance the result of Theorem~\ref{nonsmooth} by adding the
requirement that the constructed function $f$ should also be everywhere
differentiable. The construction becomes more involved, but remains explicit.
Before we proceed, we shall need the following preliminary results (lower
integral estimations for $\nu$-root type functions).

\begin{lemma}
[lower integral estimation I]\label{claim-sebastian} There exists a function
$\sigma:(0,1)\rightarrow(0,1]$ satisfying $\lim\limits_{\nu\rightarrow0}
\sigma(\nu)=0$ such that for every $x,h\in\mathbb{R}$ with $h\neq0$, we have:
\begin{equation}
\frac{1}{h}\int_{x}^{x+h}|t|^{\nu}dt\geq|x|^{\nu}\bigl(1-\sigma(\nu)\bigr).
\label{seb}
\end{equation}

\end{lemma}

\noindent\textit{Proof.} If $x=0$ the assertion follows trivially. Therefore
we may assume that $x\neq0$. Since the functions $x\rightarrow|x|^{\nu}$ are
even, we can limit our attention to the case $h>0$. We set:
\[
I=\frac{1}{h}\int_{x}^{x+h}|t|^{\nu}dt.
\]
We consider successively all four possible cases: \smallskip\newline We first assume that $x>0$ and
$t\in\lbrack x,x+h]$. In this case, $|t|^{\nu}\geq|x|^{\nu}$ and $I\geq|x|^{\nu}$,
therefore~\eqref{seb} holds for any function~$\sigma$ with nonnegative values.
\smallskip\newline Let us now assume $x<0<x+h\leq|x|$. A direct computation
gives:
\begin{equation}
I=\left(  \frac{1+y^{1+\nu}}{1+y\,\,\,\,\,\,\,}\right)  \left(  \frac
{|x|^{\nu}}{1+\nu}\right)  ,\qquad\text{where \quad}y:=\frac{x+h}{|x|}
\in\lbrack0,1]. \label{eq:bach}
\end{equation}
Consider the (continuous) functions $\Psi_{\nu}:[0,1]\rightarrow
\lbrack0,+\infty)$, $\nu\in(0,1)$, defined by
\[
\Psi_{\nu}(y)=\frac{1+y^{1+\nu}}{(1+\nu)(1+y)},\text{\quad}y\in\lbrack0,1].
\]
Then the functions $\{\Psi_{\nu}\}_{\nu>0}$ converge pointwise to the function
$\Psi\equiv1$ as $\nu$ tends to $0$. Since the above convergence is monotone,
we deduce from Dini theorem that the convergence is uniform. Setting
\begin{equation}
\sigma(\nu):=1-\min\limits_{y\in\lbrack0,1]}\frac{1+y^{1+\nu}}{(1+\nu)(1+y)}
\label{eq:bach1}
\end{equation}
we readily deduce that $\lim\limits_{\nu\rightarrow0}\sigma(\nu)=0.$ Therefore
(\ref{eq:bach}) yields
\[
I\geq|x|^{\nu}\left(  1-\sigma(\nu)\right)
\]
and (\ref{seb}) holds true for $\sigma$ given in \eqref{eq:bach1}.\smallskip\newline 
If $x<x+h\leq0$, then a direct computation yields 
$$I\,\geq\,\frac{|x|^{\nu}}{1+\nu} = \left( 1- \frac{\nu}{1+\nu}\right)|x|^{\nu}\,\geq\, (1-\sigma(\nu))|x|^{\nu},$$
where $\sigma$ is given by~\eqref{eq:bach1}.\smallskip\newline
It remains to deal with the case $x<0<|x|<x+h$. In this
case we have
\begin{align*}
I  =\frac{1}{h}\int_{x}^{|x|}|t|^{\nu}dt+\frac{1}{h}\int_{|x|}
^{x+h}|t|^{\nu}dt &\geq\,\frac{1}{h}\,\Big\{  (|x|-x)|x|^{\nu}\bigl(1-\sigma
(\nu)\bigr)+(x+h-|x|)|x|^{\nu}\Big\} \\
&  \geq\, (1-\sigma(\nu))\, |x|^{\nu}.
\end{align*}
Therefore~\eqref{seb} is still satisfied and the proof is complete. \hfill$\Box$

\bigskip

\noindent We now extend (\ref{seb}) to a more general class of functions.
Fixing parameters $d\in\mathbb{R}$, $m>0$ and $\varepsilon>0$, we set for each
$\nu\in(0,1)$
\begin{equation}
R(t)=\min\left\{  m^{\nu},\left(  \frac{|t-d|}{\varepsilon}\right)  ^{\!\nu
}\right\}  ,\qquad t\in\mathbb{R}\text{.} \label{eq:3.2}
\end{equation}
The above function is continuous and nonnegative. The following result
shows that $R$ also satisfies the same lower integral estimation as in~\eqref{seb}.

\begin{lemma}
[lower integral estimation II]\label{sebastian}For every $x,h\in\mathbb{R}$
with $h\neq0$, the function $R$ given in~\eqref{eq:3.2} satisfies
\begin{equation}
\frac{1}{h}\int_{x}^{x+h}R(t)dt\geq R(x)\bigl(1-\sigma(\nu)\bigr),
\label{eq:3.3}
\end{equation}
where $\sigma:(0,1)\rightarrow(0,1]$ is the function defined in Lemma~\ref{claim-sebastian}. \\
$($Notice that this integral estimate does not depend on
the values of the parameters $\varepsilon$, $d$ and $m.)$
\end{lemma}

\noindent\textit{Proof.} We first consider the case $\varepsilon=1$ and $d=0.$
Let $x\in\mathbb{R}$ and $h>0$. If $m\leq x<x+h$ or if $x<x+h\leq-m$ there is
nothing to prove since in both cases the function $R$ is constant on the
interval $[x,x+h]$. The case $-m\leq x<x+h\leq m$ follows from the previous
lemma, since in this case $R(t)=|t|^{\!\nu}$ on $[x,x+h]$.\smallskip\newline
Let us now consider the case $-m\leq x\leq m<x+h$. Then, according to the
previous lemma, $\int_{x}^{m}R(t)dt\geq R(x)\bigl(1-\sigma(\nu)\bigr)(m-x)$.
Since
\[
\int_{m}^{x+h}R(t)dt=m^{\nu}(x+h-m)\geq R(x)(x+h-m),
\]
we deduce
\[
\frac{1}{h}\int_{x}^{x+h}R(t)dt=\frac{1}{h}\Bigl(\int_{x}^{m}R(t)dt+\int
_{m}^{x+h}R(t)dt\Bigr)\geq R(x)\bigl(1-\sigma(\nu)\bigr).
\]
It remains to consider the case $x<-m<x+h$. In this case
\[
\int_{x}^{-m}R(t)dt=R(x)(-m-x)
\]
and
\[
\frac{1}{h}\int_{x}^{x+h}R(t)dt=\frac{1}{h}\Bigl(\int_{x}^{-m}R(t)dt+\int
_{-m}^{x+h}R(t)dt\Bigr)\geq R(x)\bigl(1-\sigma(\nu)\bigr)
\]
since, according to the previous case,
\[
\int_{-m}^{x+h}R(t)dt\geq R(-m)\bigl(1-\sigma(\nu
)\bigr)(x+h+m)=R(x)\bigl(1-\sigma(\nu)\bigr)(x+h+m).
\]
This proves the validity of~\eqref{eq:3.3} for the function $R(x)=\min\left\{
m^{\nu},|x|^{\nu}\right\}  .$ The general case for arbitrary
values of the parameters $d\in\mathbb{R}$ and $\varepsilon>0$ in~\eqref{eq:3.2} easily follows 
by translation and a standard argument.\hfill$\Box$

\bigskip

We shall also need the following refinement of Lemma~\ref{Cantor0}.

\begin{lemma}
\label{Cantor} Let $C$ be any compact subset of $[0,1]$ which is homeomorphic
to the Cantor set $\Delta:=\{0,1\}^{\mathbb{N}}$. Then there exists a
continuous curve $\gamma:[0,1]\rightarrow\mathbb{R}^{2}$, such
that
\[
\gamma\bigl([0,1]\bigr)=\gamma(C)=\mathbb{T}^{+}\qquad(\text{see
\eqref{q:T+}})
\]

\end{lemma}

\noindent\textit{Proof.} Let $\varphi$ be a homeomorphism from $C$ onto
$\{0,1\}^{\mathbb{N}}$ and let $\psi:\{0,1\}^{\mathbb{N}}\rightarrow
\lbrack0,1]$ be defined as follows:
\[
\psi\bigl((x_{n})\bigr)=\sum\limits_{n\geq1}2^{-n}x_{n}.
\]
It follows easily that $\psi$ is continuous and surjective, therefore,
$\gamma_{2}:=\psi\circ\varphi$ is a continuous function from $C$ onto~$[0,1]$.
By Urysohn lemma, we can extend $\gamma_{2}$ to a continuous curve
$\tilde{\gamma}_{2}$ from~$[0,1]$ onto~$[0,1]$. If $\gamma_{1}$ denotes the
function constructed in Lemma~\ref{Cantor0}, then the continuous curve
$\gamma:=\gamma_{1}\circ\tilde{\gamma}_{2}$ satisfies the assertion.\hfill
$\Box$\newline

\bigskip

We are now ready to construct the desired function $f$.

\begin{theorem}
[smooth exhaustive function in $\mathbb{R}$]\label{1D}There exists a
$1$-Lipschitz differentiable function $f:\mathbb{R}\rightarrow\mathbb{R}$ with
compact support for which the range of its Clarke subdifferential contains all
closed sub-intervals and all singletons of $[-1,1].$
\end{theorem}

\noindent\textit{Proof.} Let $\mathcal{D}=\{d_{n}\}_{n\geq1}$ be a countable
dense subset of $[0,1]$. Let $\{\varepsilon_{n}\}_{n}$ be a nonincreasing
sequence of positive real numbers such that $\sum\limits_{n\geq1}
\varepsilon_{n}<1/2$. Let $\{\nu_{n}\}_{n}$ be a sequence in $(0,1)$ such that
$\sum\limits_{n\geq1}\sigma(\nu_{n})<+\infty$, where $\sigma$ is the function
defined by~\eqref{eq:bach1} and evoked in Lemma~\ref{sebastian}. Let $r_{n}$ be the function defined by
\begin{equation}
r_{n}(x):=\Bigl(\frac{|x-d_{n}|}{\varepsilon_{n}}\Bigr)^{\nu_{n}},\qquad
x\in\mathbb{R}\text{.} \label{eq:r_n}
\end{equation}
According to our choice of $\varepsilon$, the set
\[
F:=[0,1]\diagdown\left(  \bigcup\limits_{d\in\mathcal{D}}\left(
d_{n}-\varepsilon_{n},d_{n}+\varepsilon_{n}\right)  \right)
\]
is a closed subset of $[0,1]$ of positive Lebesgue measure. Since $F$ is an
uncountable Polish space, there exists a closed subset $C$ of $F$ which is
homeomorphic to the Cantor set $\Delta=\{0,1\}^{\mathbb{N}}$ {(see \cite[Corollary~ 6.5]{Kechris} e.g.)} Notice also
that
\[
r_{n}(x)\geq1,\qquad\text{for all}\;x\in F\,\text{ and }\;n\in\mathbb{N}.
\]
We are now ready to construct our function $f$. Let
\[
\gamma(x)=\bigl(\alpha(x),\beta(x)\bigr)\in\mathbb{T}^{+},\qquad x\in
\lbrack0,1],
\]
be the continuous curve constructed in Lemma~\ref{Cantor} with respect to the
closed subset $C$ of $F$ evoked above. We set
\[
g_{0}(x)=\beta(x)
\]
and define inductively
\[
g_{n}(x)=\min\big\{g_{n-1}(x),\alpha(x)+r_{n}(x)\big\},\qquad\text{for }
n\geq1.
\]
Finally, we set
\begin{equation}\label{eq:dev}
g(x)=\inf\limits_{n\geq1}\,g_{n}(x)=\min\big\{\beta(x),\,\alpha(x)+\inf\limits_{n\geq1}
\,r_{n}(x)\big\} \qquad\text{and}\qquad f(x)=\int_{0}^{x}g(t)dt.
\end{equation}
Notice that the function $g$ is upper semi-continuous (as infimum of
continuous functions), hence measurable, with values in $[0,1]$ because $\alpha\leq g\leq \beta$. 
Therefore, the function $f$ is $1$-Lipschitz and nondecreasing. By construction, we have
\[
g(x)-\alpha(x)\leq r_{n}(x),\qquad\text{for every }x\in\lbrack0,1]\;\text{and
}\;n\geq1.
\]
Let us fix $x\in\lbrack0,1]$ and define
\[
R_{n}(t):=\min\{g(x)-\alpha(x),r_{n}(t)\},\qquad\text{for all }t\in
\lbrack0,1].
\]
It follows readily that $R_{n}(x)=g(x)-\alpha(x)$, thus $0\leq R_n(x)\leq \beta(x)\leq 1$. Since
\[
\max\{0,g(x)-\alpha(x)-r_{n}(t)\}=R_{n}(x)-R_{n}(t),
\]
we obtain from~\eqref{eq:3.2}--\eqref{eq:3.3} with $d=d_n$, $\nu=\nu_n$ and $m=\bigl(g(x)-\alpha
(x)\bigr)^{1/\nu_{n}}$
\begin{equation}
\frac{1}{h}\int_{x}^{x+h}\max\{0,g(x)-\alpha(x)-r_{n}(t)\}dt=R_{n}(x)-\frac
{1}{h}\int_{x}^{x+h}R_{n}(t)dt\leq R_{n}(x)\sigma(\nu_{n})\leq\sigma(\nu_{n}).
\label{r_n}
\end{equation}

\noindent\textit{Claim 1}: The function $f$ is differentiable at every point and $f^{\prime
}=g$.\medskip\newline 
\textit{Proof of Claim~1.} We shall consider separately two
cases:\medskip\newline\noindent 
--- Case $g(x)=\alpha(x)$. \smallskip\newline
Since $g\geq\alpha$, $g(x)=\alpha(x)$, $g$ is upper semi-continuous and
$\alpha$ is continuous, we deduce that $g$ is continuous at $x$, and
consequently $f$ is differentiable at $x$ with $f^{\prime}(x)=g(x)$. Notice
that the level set
\[
\lbrack g-\alpha=0]:=\{x\in\lbrack0,1]:\;g(x)=a(x)\}
\]
of the function $g-\alpha$ is dense $\mathcal{G}_{\delta}$ in $[0,1]$: indeed, it
contains the dense set $\mathcal{D}=\{d_{n}\}_{n\geq1}$ (notice that $r_{n}
(d_{n})=0$ and consequently, by~\eqref{eq:dev}, $g(d_n)=\alpha(d_n)$, for every $n\geq1$) and it is $\mathcal{G}_{\delta}$ since the
strict sublevel sets
\[
\Big[ g-\alpha<\frac{1}{n}\Big ]=\left\{  x\in\lbrack0,1]:\;g(x)-\alpha(x)<\frac{1}
{n}\right\}
\]
are open (thanks to the upper semicontinuity of $g$ and the continuity of $\alpha$)
and
\[
\Big[ g-\alpha=0 \Big]=\bigcap\limits_{n\geq1}\big[ g-\alpha<1/n\big].
\]
\smallskip\noindent --- Case $g(x)>\alpha(x)$. \smallskip\newline Since $g$ is
upper semi-continuous, we always have
\[
\limsup_{h\rightarrow0}\,\frac{f(x+h)-f(x)}{h}=\limsup_{h\rightarrow0}\,\frac
{1}{h}\int_{x}^{x+h}g(t)dt\leq g(x).
\]
It remains to prove that for fixed $\varepsilon>0$, there exists $h_{1}>0$
such that, if $|h|\leq h_{1}$, then
\begin{equation}
\frac{f(x+h)-f(x)}{h}\,\geq\, g(x)-5\varepsilon. \label{lower}
\end{equation}
Without loss of generality, we may assume that
\[
\kappa:=g(x)-\alpha(x)-\varepsilon>0.
\]
Thus, for any $n\in \mathbb{N}$, $x$ does not belong to the closed set $r_{n}^{-1}\bigl(\{\kappa\}\bigr)$,
which yields that
\[
\mathrm{dist}\bigl(x,r_{n}^{-1}\bigl(\{\kappa\}\bigr)\bigr)\,=\,\mathrm{dist}\bigl( x, r_n^{-1}([0,\kappa])\bigr)\,>\,0.
\]
Moreover, up to a subsequence,
\[
\mathrm{dist}\bigl(x,r_{n}^{-1}\bigl(\{\kappa
\}\bigr)\bigr)\underset{n\rightarrow+\infty}{\longrightarrow}0.
\]
Therefore, setting
\[
N(x,h):=\min\big\{n\geq1:\;r_{n}^{-1}\bigl(\{\kappa\}\bigr)\cap\lbrack
x-h,x+h]\neq\emptyset\big\},\qquad\text{for }h>0
\]
we deduce easily that
\[
\lim\limits_{h\rightarrow0}\,N(x,h)=+\infty.
\]
Let us fix $h_{0}>0$ such that $N:=N(x,h_{0})$ satisfies both
\[
\sum\limits_{n>N}\sigma(\nu_{n})<\varepsilon\qquad\text{and}\qquad\text{
}|g_{N}(x)-g(x)|<\varepsilon.
\]
Then, we fix $0<h_{1}\leq h_{0}$ such that, if $t\in\lbrack x-h_{1},x+h_{1}]$,
then
\[
|g_{N}(x)-g_{N}(t)|\leq\varepsilon\qquad\text{and}\qquad|\alpha(x)-\alpha
(t)|\leq\varepsilon.
\]
Consequently, if $|h|\leq h_{1}$, we have
\[
\frac{1}{h}\int_{x}^{x+h}g_{N}(t)dt\geq g_{N}(x)-\varepsilon\geq
g(x)-\varepsilon.
\]
Therefore, in order to prove~\eqref{lower}, it is enough to prove that
\[
\frac{1}{h}\int_{x}^{x+h}\bigl(g_{N}(t)-g(t)\bigr)dt\leq 4\varepsilon
,\qquad\text{whenever }|h|\leq h_{1}.
\]
Since $g(t)=\min\big\{g_{N}(t),\alpha(t)+\inf\limits_{n>N}r_{n}(t)\big\},$ we
obtain
\[
g_{N}(t)-g(t)=\max\,\left\{  \,0,\;\sup_{n>N}\{g_{N}(t)-\alpha(t)-r_{n}
(t)\}\right\}  .
\]
If $|t-x|\leq h_{1}$, then we also have that $g_{N}(t)-\alpha(t)\leq
g(x)-\alpha(x)+3\varepsilon$. Hence,
\[
\begin{split}
g_{N}(t)-g(t)  &  \leq\max\,\left\{  \,0,\;\sup_{n>N}\{g(x)-\alpha
(x)+3\varepsilon-r_{n}(t)\}\right\} \\
&  \leq\max\,\left\{  0,\;\sup_{n>N}\{g(x)-\alpha(x)-r_{n}(t)\}\right\}
+3\varepsilon\\
&  \leq\sum_{n>N}
\max\,\left\{  0,\;\{g(x)-\alpha(x)-r_{n}(t)\}\right\}  +3\varepsilon
\end{split}
\]
Integrating the above inequality, we obtain thanks to~\eqref{r_n}
\begin{equation}
\begin{split}
\frac{1}{h}\int_{x}^{x+h}\bigl(g_{N}(t)-g(t)\bigr)dt  & \, \leq \, \sum_{n>N}
\frac{1}{h}\int_{x}^{x+h}\max\left\{  0,\;g(x)-\alpha(x)-r_{n}(t)\right\}
dt+3\varepsilon\\
&  \,\leq\, \sum_{n>N}\sigma(\nu_{n})+3\varepsilon\, \leq\, 4\varepsilon
\end{split}
\end{equation}
Thus, we have shown that $f$ is differentiable at each point and that
$f^{\prime}=g$. \hfill$\lozenge$\newline

\bigskip

\noindent\textit{Claim 2}: $\operatorname{Im}(\partial f)=\partial
f({(0,1)})=\mathcal{K}_{(0,1)}:=\big\{[a,b]:\;\,0\leq a\leq b\leq1\big\}$
.\medskip\newline 
\textit{Proof of Claim~2.} For every $x\in (0,1)$, we have $0\leq f^{\prime
}(x)=g(x)\leq1$, whence $\partial f(x)\subset\lbrack0,1]$. Let us now fix
$x\in C$. Since $C\subset F$, we have $r_{n}(x)\geq1$ for all $n\geq1$ and
consequently
\[
f^{\prime}(x)=g(x)=\beta(x)\in\partial f(x).
\]
Since the set $[f^{\prime}=\alpha]=[g=\alpha]$ is dense in $[0,1]$ and
$\alpha$ is continuous, we deduce that $\alpha(x)\in\partial f(x)$, hence
$[\alpha(x),\beta(x)]\subset\partial f(x)$. 
{The reverse inclusion follows easily from~\eqref{clarke}, since $f^{\prime}(x)=g(x) \in [\alpha(x),\beta(x)]$ and the functions $\alpha$ and $\beta$ are continuous.}\smallskip\newline 
Let us finally recall that the curve $\gamma=(\alpha,\beta)$ satisfies the conclusion of
Lemma~\ref{Cantor}. This ensures that
\[
\partial f(C)=\big\{[a,b]:\;\,0\leq a\leq b\leq1\big\}.
\]
We conclude that $\partial f((0,1))=\partial f(C)=\mathcal{K}_{[0,1]}$ as asserted.\hfill$\lozenge$ \smallskip\newline
Replacing again $f$ by $\widetilde{f}:=2f-I$, where $I$ is the identity on
$[0,1],$ we obtain a differentiable function $\widetilde{f}$ with derivatives
in $[-1,1].$ It easily follows that $\widetilde{f}$ is $1$-Lipschitz and
satisfies $$\partial\widetilde{f}({(0,1)})=\mathcal{K}_{[-1,1]}
=\big\{[a,b]:\;\,-1\leq a\leq b\leq1\big\}.$$ The proof is complete.\hfill
$\Box$

\bigskip

\begin{remark} \label{dense} (i). In the above construction, $\partial f(x)$ is a singleton if and only if $x$
belongs to the ($\mathcal{G}_{\delta}$ dense) subset $[g=\alpha]$ of $[0,1]$. Moreover,
since $g=f^{\prime}$ has the Darboux property, we can easily deduce that $\partial f(x) =\partial_L f(x)$,
for all $x\in (0,1)$ and consequently, the conclusion also holds for the limiting subdifferential. \smallskip\newline

\noindent(ii). We can assume that $C$ is contained in $(0,1)$ and that $\alpha(0)=\beta(0)=\alpha(1)=\beta(1)=0$.
This allows to extend $f$ to a differentiable function on $\mathbb{R}$
satisfying $f^{\prime}(0)=f^{\prime}(1)=0$.\smallskip\newline

\noindent(iii). We can also assume that
$\mathcal{L}_{1}(C)=0$. In this case we have a negligible set $C$ satisfying
$$\partial f(C):=\{\partial f(x):x\in C\}=\mathcal{K}_{[-1,1]}.$$ \smallskip

\noindent(iv). It is clear from the above construction that the
domain of $f$ can be any nontrivial interval of arbitrarily small length and
that the range of $f$ can be taken inside $[c-\varepsilon,c+\varepsilon]$ 
for any choice of $c\in\mathbb{R}$ and $\varepsilon>0$. It follows easily, by a
standard argument, that for any nonempty open interval $\mathcal{J}$ of $\mathbb{R}$,
the set of bounded, differentiable, Clarke exhaustive $k$-Lipschitz functions in
$\mathcal{J}$ is $d_{\infty}$-dense in the (complete) metric space $\left(\mathrm{Lip}^{[k]}(\mathcal{J}),d_{\infty}\right)$ 
of all bounded Lipschitz functions in~$\mathcal{J}$ with $\|f\|_{\mathrm{Lip}}\leq k$. \smallskip\newline

\noindent(v). Let $I=(a,b)$ be a nonempty (possibly unbounded) interval. Then the set $\mathcal{E}$ of all real-valued Lipschitz functions in $I$ which are everywhere differentiable and Clarke
exhaustive is spaceable when equipped with the semidistance $d_{\mathrm{Lip}}(f,g):=\| f-g \|_{\mathrm{Lip}}$, for all $f$, $g$ in $\mathcal{E}$. \smallskip\newline
Indeed, it is sufficient to consider a sequence of disjoint intervals $\{(a_n,b_n)\}_n$ such that 
\[
a<a_n<b_n<a_{n+1}<b\,, \quad \text{for every } n\in \mathbb{N}\,, 
\] a sequence of Clarke exhaustive functions $\{f_n\}_n$ such that $\|f_n\|_{\mathrm{Lip}}=1$ and $\mathrm{supp}\,f_n\subset (a_n,b_n)$ for all $n\in \mathbb{N}$, and the operator
$T:c_0(\mathbb{N})\to \mathcal{E}$ defined by
\[T(\{x_n\}_n):= \sum_{n=1}^\infty x_nf_n(\cdot).\]
Since the supports of the functions $f_n$ are pairwise disjoint, it follows easily that the operator $T$ is well defined and establishes a linear isometry between $c_0(\mathbb{N})$ and its image. 
Therefore, the metric space $(\mathcal{E},d_{\mathrm{Lip}})$ contains an isometric copy of $c_{0}(\mathbb{N})$. Similar constructions of operators $T$ can be found in \cite{DF2019, BFT2021}.
\smallskip \newline

\noindent (vi). The set of all Lipschitz functions in $[0,1]$
which are everywhere differentiable and Clarke exhaustive cannot be $\Vert
\cdot \Vert _{\infty }$-spaceable in $\left( \mathrm{Lip}([0,1]),\Vert \cdot
\Vert _{\mathrm{\infty }}\right) $ (the latter being seen as a dense
subspace of the Banach space $(\mathcal{C}([0,1],\Vert \cdot \Vert _{\mathrm{\infty }})$). This is a straightforward consequence of the classical fact
that every subspace $Y$ of Lipschitz functions which is $\Vert \cdot \Vert _{\mathrm{\infty }}$-closed in $\mathcal{C}([0,1])$ is necessarily finite
dimensional. Let us sketch a proof for reader's convenience: we consider the family of linear operators $\{T_{x,y}:\,x,y\in \lbrack 0,1],\,x\neq y\}$ 
defined by $T_{x,y}(f)=\frac{f(x)-f(y)}{|x-y|},$ for all $f\in Y\subset \mathrm{Lip}([0,1])$. Since $T_{x,y}(f)\leq ||f||_{\mathrm{Lip}}$ for all $x,y\in \lbrack 0,1],\,x\neq y,$
and $(Y,\Vert \cdot \Vert _{\mathrm{\infty }})$ is complete, applying the Banach-Steinhaus theorem we deduce that for some $M>0$ and all $x,y\in
\lbrack 0,1],\,x\neq y$, it holds $||T_{x,y}||\leq M.$ It follows from Arzel\`{a}-Ascoli theorem that every $\vert|\cdot \vert|_{\infty }$-bounded sequence $\{f_{n}\}_{n}$ in $Y$ has a converging subsequence, and consequently, the closed unit ball $\overline{B}_{Y}(0,1)$ of $Y$ is compact, ensuring that $Y$ is
finite dimensional.
\end{remark}


\subsection{Subdifferential containing all convex bodies in
$\mathbb{R}^{N}$ ($N\geq 2$)} \label{ss-3.2}

We shall now deal with the higher dimensional case and construct a differentiable Lipschitz
function $f$ which is almost exhaustive, that is, its Clarke subdifferential
contains all nonempty convex compact subsets of $\overline{B}(0,\Vert
f\Vert_{\mathrm{Lip}})$ of nonempty interior. The question of whether it is
possible to obtain a Lipschitz Clarke exhaustive function in dimension $N\geq 2$
remains open.\smallskip
 
Let us stress the fact that the forthcoming
construction cannot be applied in one dimension. Roughly speaking, our approach occupies
one-dimension to code the family of convex bodies in $\overline{B}(0,\Vert f\Vert_{\mathrm{Lip}})$  
(based on the fact that any compact geodesic metric space can be represented as a continuous surjective 
image of $[0,1]$) and requires at least one extra dimension to make an efficient use of this coding.
Although the overall construction is less explicit and more involved, the reader can possibly trace some analogies 
between the aforementioned surjection and the curve obtained in Lemma~\ref{Cantor0} which was used to recover all
closed intervals in $[0,1]$.\smallskip
 
In order to keep notation simple, $\mathbb{R}^{N}$ will be considered with its
natural Euclidean structure (despite the fact that our results
Lemma~\ref{lem: function Rn} and Theorem~\ref{theo: function RN} hold true in
any finite dimensional normed space). Therefore, by Riesz representation
theorem, the dual space of $\mathbb{R}^{N}$ will be identified to itself. We shall also identify
1-forms $Df(x)$ with gradients $\nabla f(x)$, for any differentiable function $f:\mathbb{R}^N\rightarrow\mathbb{R}$. 
In what follows we are going to construct:

\begin{itemize}
\item[(I)] for every $n\geq 1$, a compactly supported differentiable $n$-Lipschitz function
$f_{n}:\mathbb{R}^{N}\rightarrow\mathbb{R}$ whose Clarke subdifferential
contains in its range every compact convex subset of nonempty interior that
lie in the closed ball $\overline{B}(0,n).$ 
\end{itemize}

\noindent Similarly to the one-dimensional case, the method of construction will
directly yield that the set of all functions as above is $d_{\infty}$-dense in
$\left(  \mathrm{Lip}(\mathcal{U}),\|\cdot\|_{\infty}\right)  $ (for
$\mathcal{U}\subset\mathbb{R}^{N}$ open and bounded) and $\|\cdot
\|_{\mathrm{L}}$-spaceable in $\left(  \mathrm{Lip}(\mathcal{U}),\|\cdot
\|_{\mathrm{L}}\right).$

\begin{itemize}
\item[(II)] a differentiable locally Lipschitz function $f:\mathbb{R}
^{N}\rightarrow\mathbb{R}$ whose Clarke subdifferential contains in its range
all compact convex bodies of $\mathbb{R}^{N}$.
\end{itemize}

\medskip

\noindent Notice that the second assertion follows directly from the first: it is enough
to consider a family of differentiable Lipschitz functions $f_{n}:\mathbb{R}
^{N}\rightarrow\mathbb{R}$ with $\mathrm{Lip}(f_n)=n$ and disjoint supports (for instance,
$\mathrm{supp}(f_{n})\subset B(3n\,e_{1},1)$ where $e_{1}=(1,0,\cdots,0)$),
satisfying the statement (I) and define the function
\begin{equation}\label{eq:AD}
f(x)={\displaystyle\sum\limits_{n\geq1}}f_{n}(x),\quad\text{for all }
x\in\mathbb{R}^{N}.
\end{equation}
One readily gets that $f$ is everywhere differentiable, locally Lipschitz and satisfies assertion~(II).
\smallskip

\noindent Let us now proceed to the construction evoked in (I). It clearly suffices to
do it for the case $n=1$ and construct a $1$-Lipschitz function. \\ 
This will be done in two stages: we first fix a compact convex subset $C$ in $\mathbb{R}^{N}$ that contains $0$ and construct an $L$-Lipschitz function (with $C\subset \overline{B}(0,L)$)
whose Clarke subdifferential contains all compact convex subsets $K$ of $C$
that contain $0$. The general case will follow using separability arguments, by considering an adequate sequence $\{C_{n}\}_{n}$ of
compact convex sets with $0\in\mathrm{int}\,C_{n}$, then gluing adequate
translations of the corresponding constructed functions.

\subsubsection{An intermediate construction}

For a nonempty compact convex subset $C$ of $\mathbb{R}^{N}$ with $0\in C$, recalling from~\eqref{eq:G}--\eqref{eq:DH} 
the definition of $(\mathcal{K}_{C},D_{\mathrm{H}})$ we denote by
\begin{equation}
\mathcal{K}_{C}^{0}:=\{K\in \mathcal{K}_{C}:~0\in K\}  \label{eq:Go}
\end{equation}
the set of all convex compact subsets of $C$ containing $0$. Notice that $\mathcal{K}_{C}^{0}$ is closed in $\mathcal{K}_{C}$, therefore 
$(\mathcal{K}_{C}^{0},D_{\mathrm{H}})$ is a compact metric space. Moreover, it  is a \textit{geodesic space} {(see \cite[p.~72]{S} \textit{e.g.})}. 
Indeed, for any two elements ${K_{0},K_{1}}\in\mathcal{K}_{C}^{0}$ and $\lambda\in(0,1)$, we have:
\begin{equation}
{K_{\lambda}:=(1-\lambda)K_{0}+\lambda K_{1}}\in\mathcal{K}_{C}^{0}
\qquad\text{and}\qquad D_{\mathrm{H}}({K_{0}, K_{\lambda})=\lambda D_{\mathrm{H}}
(K_{0},K_{1})}.\label{eq:veliov}
\end{equation}

\noindent We shall show, as an application of the next lemma, that there exists a differentiable $1$-Lipschitz function $f:\mathbb{R}^N\to\mathbb{R}$ such that 
$\mathcal{K}_{\overline{B}(0,1)}^0$ is contained in the image of the subdifferential of~$f$.

\begin{lemma} \label{lem: function Rn} 
Let $C\subset\mathbb{R}^{N}$ be a convex compact set such that $0\in C$ and $L:=\underset{x\in C}{\max}\{\Vert x \Vert\}.$ Then:\smallskip\newline
(i). There is a differentiable $L$-Lipschitz continuous and compactly supported function $f:\mathbb{R}^{N}\rightarrow\mathbb{R}$ such that:
\begin{equation} \label{eq:salas1}
\text{for every } K\in\mathcal{K}_{C}^{0}, \text{ there exists } x \in\mathbb{R}^{N} \text{ such that } \partial f(x)=K.
\end{equation}

\noindent (ii). Let us further assume that $0\in \mathrm{int}(C)$. Then in addition to \eqref{eq:salas1} we get: 
\begin{equation}\label{eq:salas2}
\partial f(x)\subset C, \text{ for all } x\in\mathbb{R}^{N}
\end{equation}
\end{lemma}

\begin{figure}[ht]
\centering
\includegraphics[width=10cm]{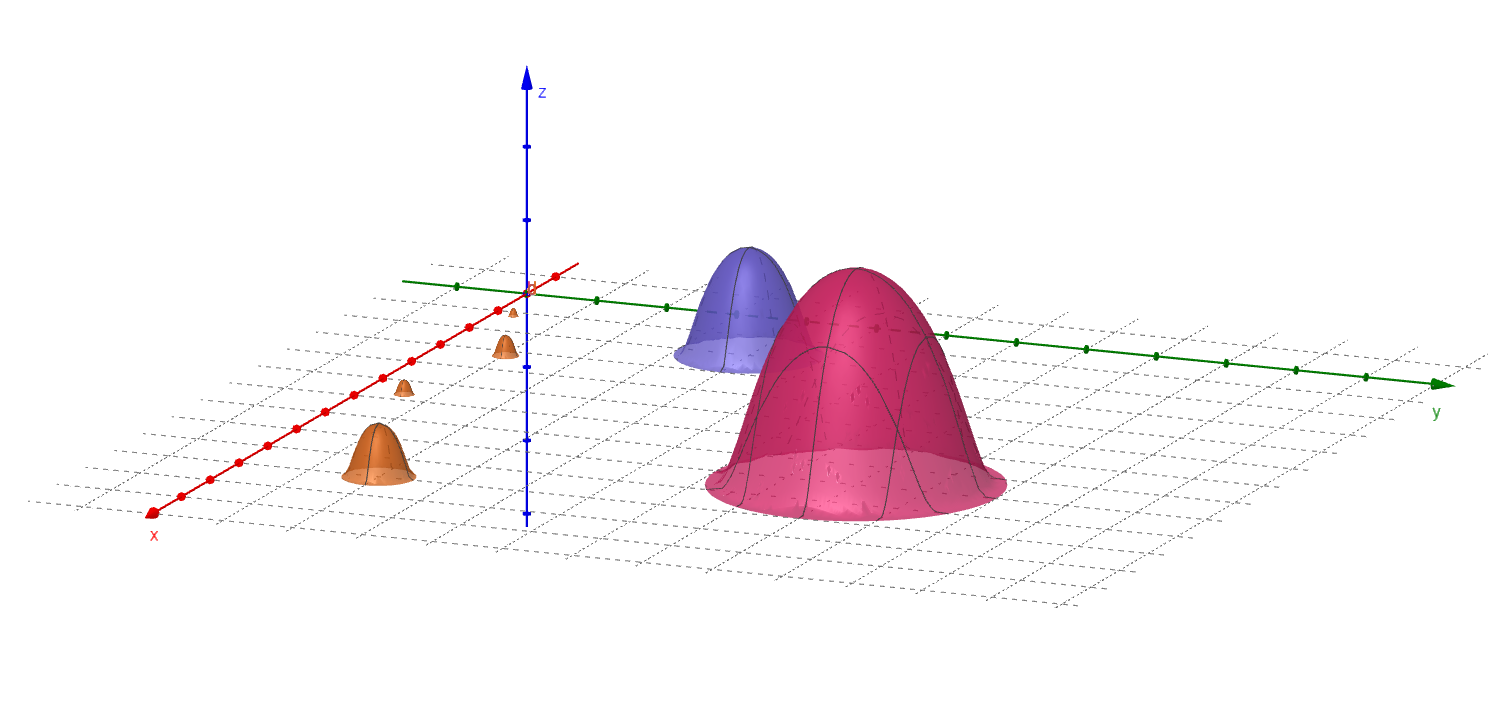}
\caption{Sketch of the function constructed in Lemma~\ref{lem: function Rn}.}
\centering
\end{figure}

\noindent\textbf{Proof.} 
\textbf{(i)}.  If $C=\{0\}$, then the function $f\equiv 0$ satisfies trivially the conclusion. Therefore, we may assume $\{0\} \subsetneq C$. 
Since $\left(\mathcal{K}_{C}^0,D_{\mathrm{H}}\right)$ is a compact metric space, there exists a continuous surjective map from the Cantor set 
$\Delta$ to~$\mathcal{K}_{C}^0$ {(see \cite[Theorem~4.18]{Kechris})}. Since $\mathcal{K}_{C}^0$ is also geodesic, a standard argument shows that this map can be extended to a continuous surjective map
$$h:[0,1]\rightarrow\mathcal{K}_{C}^0 \qquad \text{(coding the elements of $\mathcal{K}_{C}^0$)}.$$ 
Let
\[
\mathcal{D}=\left\{  d_{n}:\text{ }n\in\mathbb{N}\right\}
\]
be a countable dense subset of $(0,1)$. \smallskip\newline
Let us define, by induction, two sequences $\{\alpha_n\}_{n}$ and  $\{\varepsilon_n\}_{n}$,
satisfying $\alpha_n>\varepsilon_n>0$, for all $n\geq 1$, satisfying $$\lim\limits_{n\to\infty}\alpha_n=\lim\limits_{n\to\infty}\varepsilon_n=0$$ 
and the following property: setting
\begin{equation}
Q_{n}=\left(d_{n},\alpha_n,0,\ldots,0\right)\in\mathbb{R}^{N},
\quad\text{where } n\in\mathbb{N},\label{eq:Qin}
\end{equation}
the balls $\overline{B}(Q_{n},\varepsilon_n)$ are pairwise disjoint and contained in $(0,1)^{N}$.\smallskip\newline
Indeed, assuming that $\alpha_1,\cdots,\alpha_n,\varepsilon_1,\cdots,\varepsilon_n$ have already been constructed accordingly,
pick $$0<\alpha_{n+1}<m_n:=\min\{\alpha_i-\varepsilon_i;\,1\le i\le n\}$$ 
and then choose $0<\varepsilon_{n+1}<\alpha_{n+1}$ such that 
$\alpha_{n+1}+\varepsilon_{n+1}<m_n$ and $\varepsilon_{n+1}<\min\{d_{n+1},1-d_{n+1}\}$.
Notice that we can also assume the extra condition $\lim\limits_{n\to\infty}\varepsilon_n/\alpha_n=0$ (which will be needed later). 
{A concrete choice of such sequences is given by $\alpha_{n}=1/2^{n}$, $\varepsilon_{n}=1/n2^{n+2}$ 
and $\{d_{n}\}_{n\in\mathbb{N}}$ be a standard enumeration of the dyadics in $(0,1)$ given by $d_{1}=1/2$
and
\[
d_{n}=\frac{2\,i(n)-1}{2^{m(n)+1}},
\]
where for every $n\geq 2$ we denote by $m(n)$ the unique $m\in \{1,\ldots,n\}$ such that
\[
s_{m}:={\displaystyle\sum\limits_{k=0}^{m-1}} 2^{k}<n\leq s_{m+1}:=
{\displaystyle\sum\limits_{k=0}^{m}}
2^{k}
\]
and we set $i(n):=n-s_{m(n)}\in\{1,\ldots,2^{m(n)-1}\}$.
}
Further, for every $n\in\mathbb{N}$ we define 
\begin{equation}\label{eq:salas3}
H_n:= \left( \, h(d_{n})+\overline{B}(0,\gamma_n)\, \right)\, \bigcap \, \overline{B}(0,L),
\end{equation}
where $\{\gamma_n\}_n$ is an arbitrary sequence of positive numbers converging to $0$.
Therefore, for every $n\in\mathbb{N}$ we have 
\begin{equation}\label{eq:set inclusion}
B\left(0,\min\left\{\gamma_n,L\right\} \right)\subset H_n\subset h(d_n)+\overline{B}(0,\gamma_n).    
\end{equation}
Since $H_n$ is a convex compact subset of $\mathbb{R}^N$ such that $0\in \mathrm{int}(C)$, 
according to a consequence of a result of J.~Borwein, M.~Fabian, I.~Kortezov and P.~Loewen~\cite[Theorem~12]{BFKL} (see also T.~Gaspari~\cite{G}),
for every $n\in\mathbb{N}$, there exists a $\mathcal{C}^{1}$-smooth function $b_n:\mathbb{R}^N\to\mathbb{R}$, with support in the unit ball, 
such that $\nabla b_n(\mathbb{R}^N)=H_n$ and $\|b_n\|_{\infty}\leq 1$.
We set
\begin{equation}
\phi_{n}(x):=\varepsilon_n\cdot b_n\Bigl(\frac{x-Q_{n}}{\varepsilon_n}\Bigr)\label{eq:phiin}
\end{equation}
and observe that $\phi_{n}$ is $L$-Lipschitz and satisfies:
\[
\Vert\phi_{n}\Vert_\infty\le\varepsilon_{n}\qquad\text{ and }\qquad
\mathrm{supp}(\phi_{n})\subset\overline{B}(Q_{n},\varepsilon_n).
\]
It follows that the elements of the family $\mathcal{S}=\{\mathrm{supp}
(\phi_{n}):~n\in\mathbb{N}\}$ are pairwise disjoint and
contained in $[0,1]^{N}$. Moreover, for any $x\in\mathbb{R}^{N}$ and
$n\in\mathbb{N}$, we have
\begin{equation}
\nabla\phi_n\bigl(\overline{B}(Q_n,\varepsilon_n)\bigr)=H_n.\label{eq:salas0}
\end{equation}
Notice further that if $x\notin\mathbb{R}\times\{0\}^{N-1},$ then
$B(x,\delta)$ intersects at most one element of the family~$\mathcal{S}$ for
$\delta>0$ sufficiently small.

\noindent We are ready to define the function $f$ that satisfies our
assertion:
\begin{equation}\label{eq:def-sum}
\left\{
\begin{array}
[c]{l}
f:\mathbb{R}^{N}\rightarrow\mathbb{R}\smallskip\\
f(x)=\sum\limits_{n=1}^\infty  \phi_{n}(x) .
\end{array}
\right.
\end{equation}
Since $\mathrm{supp}(f)\subset\lbrack0,1]^{N},$ the function $f$ is
compactly supported. It follows easily that $f$ is $L$-Lipschitz and coincides
with~$\phi_{n}$ in a neighborhood of~$Q_{n}$. Therefore, $\partial
f(x)\subset\overline{B}(0,L)$, for all $x\in\mathbb{R}^{N}$. 
Moreover, since $\Vert\phi_{n}\Vert_\infty\underset{n\rightarrow\infty}{\longrightarrow}0$,
$f$ vanishes and is continuous on $\mathbb{R}\times\{0\}^{N-1}$. The next claim
yields directly \eqref{eq:salas1}.\medskip\newline
\textit{Claim 1}: For every $K\in\mathcal{K}_{C}^{0}$ there exists $x\in\lbrack0,1]\times\{0\}^{N-1}$ with
$\partial f(x)=K.$ \medskip\newline
\textit{Proof of the Claim 1}. Fix $\widehat{K}\in\mathcal{K}_{C}^{0}$ and pick any $\widehat{t}\in [0,1]$
such that $h(\widehat{t})=\widehat{K}$. Set 
$$\hat{x}=(\hat{t},0,\ldots, 0)\,\in\,\lbrack0,1]\times\{0\}^{N-1}.$$
We first show that
$\widehat{K}\subset\partial f(\hat{x}).$ Indeed, by continuity of the
function~$h$ we have
\[
\underset{t\rightarrow\widehat{t}}{\lim}\,D_{H}(h(t),\widehat{K})=0.
\]
Take a sequence $d_{k(n)}\in\mathcal{D}$ 
converging to $\widehat{t}$ 
so that $\widehat{x}=\lim_{n\rightarrow\infty} Q_{k(n)}$.  
Recalling \eqref{eq:set inclusion} we deduce that:
\begin{equation}
\lim_{n\rightarrow\infty}H_{k(n)}=\widehat{K}.\label{eq:lim}
\end{equation}
Thus, if $p\in\widehat{K}$, there exist points $x_n\in\overline{B}(Q_{k(n)},\varepsilon_n)$, $n\geq 1$, such that the sequence $\bigl\{\nabla f(x_n)\bigr\}_n$
converges to $p$. Since $\widehat{x}=\lim_{n\rightarrow\infty}
x_n$, we obtain $p\in\partial f(\widehat{x})$. This proves that
$\widehat{K}\subset\partial f(\widehat{x})$.
\smallskip\newline\noindent 
Let us now prove $\partial f(\hat{x})\subset\widehat{K}$. 
Fix $\varepsilon>0$. Since $h$ is continuous, there exists $\delta>0$ such that
\begin{equation}
h(t)\subset\widehat{K}+B(0,\varepsilon/2),\qquad\text{for all \ }
t\in(\widehat{t}-\delta,\widehat{t}+\delta)\cap\lbrack0,1].\label{eq:salas7}
\end{equation}
For $\rho>0$ sufficiently small (the exact value of $\rho$ will be fixed
later) we set:
\[
\mathcal{U}_{\rho}:=\Big[ (\widehat{t}-\frac{\delta}{2},\widehat{t}
+\frac{\delta}{2})  \times(-\rho,\rho)^{N-1}\Big]\setminus\Big[ \{0\}\times
\mathbb{R}^{N-1}\Big].
\]
Since for every $x\in\mathcal{U}_{\rho}$ there is at most one $n\in\mathbb{N}$ such that $x\in\mathrm{supp}(\phi_{n}),$ it follows
that either $\nabla f(x)=0$ (if $x$ does not belong to any element
of the family $\mathcal{S}$) or in view of
(\ref{eq:salas0})
\[
\nabla f(x)=\nabla \phi_{n}(x)\in H_{n},
\]
In this latter case, since
$d_{n}\in(\widehat{t}-\delta,\widehat{t}+\delta)$ it follows from
(\ref{eq:salas7}) and (\ref{eq:salas3}) that
\[
\nabla f(x)\in H_{n}\subset h(d_{n})+\overline{B}(0,\gamma_n)\subset\left(  \widehat{K}+B(0,\frac{\varepsilon}{2})\right) +\overline{B}(0,\gamma_n).
\]
We can take $\rho>0$ sufficiently small to ensure that $\gamma_n<\varepsilon/2$, whenever $\mathrm{supp}(\phi_{n})\cap\mathcal{U}_{\rho}\neq
\emptyset$.  Choosing $\rho>0$ in this way, we infer that
\[
\partial f(x)\subset\widehat{K}+B(0,\varepsilon),~\text{for all }
x\in\mathcal{U}_{\rho}.
\]
Since the set $\mathcal{N}:=\{0\}\times\mathbb{R}^{N-1}$ is negligible for the
Lebesgue measure, we deduce easily from the formula \eqref{clarke} of the
Clarke subdifferential that
\[
\partial f(\widehat{t}\times\{0\}^{N-1})\subset\widehat{K}+B(0,\varepsilon).
\]
Since $\varepsilon>0$ can be chosen arbitrary small, we obtain the desired
conclusion.
\medskip\noindent

\noindent\textit{Claim 2}: The function $f$ is differentiable on $\mathbb{R}^N$ \medskip\newline
\textit{Proof of the Claim 2}. Since the compact sets $\text{supp}(\phi_n)$ are disjoint subsets $\mathbb{R}^N$ 
and do not intersect the closed subset $[0,1]\times\{0\}^{N-1}$ of $\mathbb{R}^N$,
the function $f$ is $\mathcal{C}^{1}$-smooth on $\mathbb{R ^N}\setminus \left([0,1]\times\{0\}^{N-1}\right)$. 
Let us now treat the case where $x\in[0,1]\times\{0\}^{N-1}$. In this case, $f(x)=0$. Take any $y\in\mathbb{R}^N$. 
If the point $y$ does not belong to $\text{supp}(\phi_n)$ for any $n$, then $f(y)=0$, while if $y\in\text{supp}(\phi_n)$ for some $n\in\mathbb{N}$,
then we deduce from~\eqref{eq:phiin} that $\vert f(y)-f(x)\vert=\vert f(y)\vert\le\varepsilon_n\ll \Vert y-x\Vert$ because $\Vert y-x\Vert\ge\alpha_n-\varepsilon_n$
and $\lim\limits_{n\to\infty}\varepsilon_n/\alpha_n=0$. Since $\text{supp}(\phi_n)$ is compactly contained in $(0,1)^\mathbb{N}$, we conclude that $f$ is differentiable at $x$ and $\nabla f(x)=0$. \hfill$\lozenge$  \smallskip\newline
\noindent This completes the proof of (i).\bigskip\newline\noindent
\textbf{(ii).} We now assume that there exists $\lambda>0$ such that
$B(0,\lambda)\subset C$. To construct a function $f$ that satisfies
(\ref{eq:salas1})--(\ref{eq:salas2}), we replace the definition of $H_{n}$  in (\ref{eq:salas3}) by
\[
H_{n}:=\left(  h(d_{n})+\overline{B}(0,\gamma_n)\right)  \cap
C,
\]
and we proceed as before. It follows easily that $\partial f(x)\subset
C\subset B(0,L),$ for all $x\in\mathbb{R}^{N}$ (in particular $f$ is
$L$-Lipschitz) and (\ref{eq:salas1}) follows as in (i).\hfill$\Box$
\bigskip
\begin{remark}\label{rem-gap}
\textbf{(i)}. A more elementary (and self-contained) proof of Lemma~\ref{lem: function Rn} can be provided if the assumption of differentiability of $f$ is dropped.
Indeed, following the lines of the above proof, once defined the sets $H_n$ in \eqref{eq:salas3}, we can consider the functions $\widehat{b}_n:\mathbb{R}^N\to\mathbb{R}$ defined by
\[
\widehat{b}_n(x):=\min\{0, \sup_{p\in H_n}\langle p,x\rangle-c_n\},
\]
where $c_n>0$ is chosen such that $\textup{supp}(\widehat{b}_n)\subset \overline{B}(0,1)$ and $\|\widehat{b}_n\|_\infty\leq 1$. Note that, in a neighborhood of $0$, 
the function $\widehat{b}_n+c_n$ is the support function of $H_n$. The function $\widehat{b}_n$ is non-differentiable and 
$\partial \widehat{b}_n(x)\subset H_n=\partial \widehat{b}_n(0)$ for all $x\in X$. Then we define $\phi_n$ as in~\eqref{eq:phiin}, using the functions $\{\widehat{b}_n:n\in\mathbb{N}\}$ (instead of $b_n$) and the function 
$\widehat{f}:\mathbb{R}^N\to\mathbb{R}$ as in~\eqref{eq:def-sum}. Proceeding as in the above proof and using the fact that the Clarke subdifferential $\partial \widehat{f}$ is outer semicontinuous and $\partial \widehat{f}(Q_n)=H_n$ for all $n\in \mathbb{N}$, we deduce that $\partial \widehat{f}((t,0,...0))=h(t)$, for any $t\in [0,1]$. \medskip\newline
\noindent \textbf{(ii)}. We can also use this idea to construct an everywhere differentiable function satisfying Lemma~\ref{lem: function Rn}. 
 Indeed, fixing a positive mollifier $\rho:\mathbb{R}^N\rightarrow\mathbb{R}$, we set $\rho_n(\cdot):=\nu_n^{-N}\rho(\cdot/\nu_n)$,  $\nu_n \geq 1$ and consider the convolution 
$b_n^\ast:= \widehat{b}_n\ast \rho_{n}$. Then taking $\nu_n>0$ sufficiently small, we ensure that $b^\ast_n$ is a good approximation of $\widehat{b}_n$, which becomes better and better as $\nu \to 0$. (The interested reader is invited to work out the details of this construction.) \medskip\newline
{\noindent\textbf{(iii)}. A careful inspection of the proof of Lemma~\ref{lem: function Rn} reveals that one can work directly with the continuous surjective map $h:\Delta\mapsto\mathcal{K}^0_C$ by simply replacing $[0,1]$ by $\Delta$ in the proof and by taking a countable dense subset $\mathcal{D}$ of 
$\Delta\setminus \{0,1\}$. The coding over $\Delta$ does not use the fact that the space $(\mathcal{K}^0_C, D_{\mathrm{H}})$ is a geodesic space. This remark will be particularly relevant in Section~\ref{ss-tapia}.}
\end{remark}

\subsubsection{Main result: recovering convex bodies}

Based on Lemma~\ref{lem: function Rn} (which recovers all convex bodies containing $0$), we can now deduce the general case. We shall also need the following lemma. 

\begin{lemma} \label{bump}  Let $x^*\in\mathbb{R}^N$ be such that $\Vert x^*\Vert<1$. Then, 
there exists a continuously differentiable and $1$-Lipschitz function $h:\mathbb{R}^N\to\mathbb{R}$ 
with support in the unit ball $B(0,1)$ and $\delta >0$ such that 
\[
\nabla h(x)=x^*(x),\qquad \text{ for all } \,x\in B(0,\delta).
\]
\end{lemma}

{\noindent The proof of the above lemma is straightforward. It is sufficient to set $x\mapsto x^{\ast }(x)$ on a small ball centered at $0$, consider an affine interpolation outside this ball which brings to the value to $0$, and finally use a mollifier with a sufficiently small support.}


\bigskip

\noindent We are now ready to state the main result of this section.

\begin{theorem}
[almost exhaustive function in $\mathbb{R}^{N}$]
\label{theo: function RN} There exists a differentiable $1$-Lipschitz compactly supported
function $f:\mathbb{R}^{N}\rightarrow\mathbb{R}$ such that for every convex body $K$ of $\overline{B}(0,1)$, there exists
$x\in\mathbb{R}^{N}$ such that $\partial f(x)=K$.
\end{theorem}

\noindent\textbf{Proof}.
Let $\{q_{n}^{\ast}\}_{n}\subset\,B(0,1)$ be a dense sequence in $\overline{B}(0,1)$. 
We claim that there exists a differentiable, $1$-Lipschitz and compactly supported function 
$g:\mathbb{R}^{N}\rightarrow\mathbb{R}$ satisfying that for each $n\in\mathbb{N}$, there
exists a set $\mathcal{U}_{n}\subset\mathbb{R}^{N}$ with nonempty interior, such that
\[
\nabla g(x) = q_{n}^{\ast}\qquad\text{if \ }x\in\mathcal{U}_{n}.
\]
Let us present a quick construction of the function $g$. First, applying Lemma~\ref{bump},
for any $n\in\mathbb{N}$,  there exists a continuously differentiable and $1$-Lipschitz function $g_n:\mathbb{R}^N\to\mathbb{R}$  with support in the unit ball such that $\nabla g_n(x)=q_n^*(x)$, for all $x$ in a neighborhood of $0$. 
Take any sequence $\{x_n\}_n$ of distinct points of $B(0,1)$ that converges to some point $\ell$ of the open unit ball, with $\ell \neq x_n$, for all $n\in \mathbb{N}$. 
Choose further $\{\varepsilon_n\}_n\subset (0,1)$ such that $\bigl\{{\overline{B}}(x_n,\varepsilon_n)\bigr\}_n$ is a sequence of disjoint closed balls contained in~$B(0,1)$. The required function $g$ is defined by
\[
g(x):=\sum_n\varepsilon_n g_n\Bigl(\frac{x-x_n}{\varepsilon_n}\Bigr)
\]
The function $g$ is the sum of disjointly supported functions, hence $g$ is $1$-Lipschitz, 
the support of $g$ is contained in the unit ball, $\nabla g(x)=q_n^*$ in a neighborhood $\mathcal{U}_{n}$ of $x_n$, 
and $g$ is differentiable at every point $x\in\mathbb{R}^N \setminus\{\ell\}$.  \smallskip\newline
Let us now show that the function $g$ is also differentiable at $\ell$, provided the sequence $\{\varepsilon_n\}_n$ satisfies
$$\lim_{n\to\infty}\frac{\varepsilon_n}{\Vert x_n-\ell\Vert}=0 \qquad\text{(by shrinking the values of $\varepsilon_n$ we can always guarantee this.)}$$
Indeed, for $n$ sufficiently large and for any $x\in B(x_n,\varepsilon_n)$, we have
\[
g(x)=\varepsilon_n\, g_n\Bigl(\frac{x-x_n}{\varepsilon_n}\Bigr) \le \varepsilon_n
\quad\text{and}\quad \frac{g(x)-g(\ell)}{\|x-\ell \|} \leq \frac{\varepsilon_{n}}{ \| x-\ell \|} \leq  
\frac{\varepsilon_{n}}{ \| x_n-\ell \|} \underbrace{\left( \frac{\| x_n -\ell \|}{ \| x-\ell \|} \right)}_{\ge 1/2}
 \underset{n \to \infty}{\longrightarrow}\, 0\,,
\]
yielding that $g$ is differentiable at $\ell$ with $\nabla g(\ell)=0$. \medskip\newline
For each $n\in\mathbb{N}$, let $x_{n}\in\mathbb{R}^{N}$ and $\lambda_{n}>0$ be
such that $\overline{B}(x_{n},\lambda_{n})\subset\,$\textup{$\mathrm{int\,}$
}$\mathcal{U}_{n}$. Set $C_{n}:=\overline{B}(-q_{n}^{\ast},1)$ and notice that
$0\in$\,\,\textup{$\mathrm{ int}\,$}$C_{n}$. Applying Lemma~\ref{lem: function Rn}(ii) for $C=C_n$,
we obtain a differentiable 1-Lipschitz function $f_n:\mathbb{R}^N\to\mathbb{R}$ satisfying 
(\ref{eq:salas1})--(\ref{eq:salas2}).
Up to a suitable
re-scaling, namely replacing $f_{n}$ by $\delta_n f_{n}(\cdot/\delta_n)$, we can
assume that \textup{$\mathrm{supp}$}$~f_{n}\subset$\,\,\textup{$\mathrm{ int}$
}$B(0,\lambda_{n})$. We define the function
\[
\left\{
\begin{array}
[c]{l}%
f:\mathbb{R}^{N}\rightarrow\mathbb{R\medskip}\\
f(x)=g(x)+\sum\limits_{n=1}^{\infty}f_{n}(x-x_{n}).
\end{array}
\right.
\]
Notice that, for any $x\in\mathbb{R}^{N}$, there is at most one $n\in
\mathbb{N}$ such that $x-x_{n}\in\,\,\,$\textup{$\mathrm{supp}$}$~{f_{n}}$.
Moreover, for any $n\in\mathbb{N}$ and $x\in\mathcal{U}_{n}$, we deduce that
\[
\nabla f(x)=q_{n}^{\ast}+\nabla f_{n}(x-x_{n})\in B(0,1),
\]
and if $x$ is not in any $\mathcal{U}_{n}$, then $\nabla f(x)=\nabla g(x)\in B(0,1)$.
It follows easily that $f$ is $1$-Lipschitz. Let us now verify that $f$
satisfies the property asserted in the statement of the theorem. To this end,
let $K\subset B(0,1)$ be a convex compact set with nonempty interior. Since
$\{q_{n}^{\ast}\}_{n}$ is dense in~$B(0,1)$, there exists $n\in\mathbb{N}$ such
that $q_{n}^{\ast}\in\mathrm{int}K$. Therefore, $K-q_{n}^{\ast}\subset
\overline{B}(-q_{n}^{\ast},1)=C_n$. From property (\ref{eq:salas1}) of
Lemma~\ref{lem: function Rn}, there exists $y\in \mathrm{supp}(f_{n})\subset\,B(0,\lambda_{n})$ 
such that $\partial f_{n}(y)=K-q_{n}^{\ast}$. Recalling that $x_{n}\in\mathcal{U}_{n}$, setting $x_{K}:=y+x_{n}\in
B(x_{n},\lambda_{n})\subset\mathcal{U}_{n}$, we obtain
\[
\partial f(x_{K})=\nabla g(x_{K})+\partial f_{n}(y)=K.
\]
The proof is complete. \hfill$\Box$

\medskip 

\begin{remark}\label{devil}
A careful inspection of the proof of Theorem~\ref{theo: function RN} reveals that for the constructed function $f:\mathbb{R}^N\rightarrow \mathbb{R}$, 
the Clarke subdifferential $\partial f(x)$ and the limiting subdifferential $\partial_L f(x)$ coincide at every point. Let us recall that the same 
situation occurred in Theorem~\ref{1D}, for the case $N=1$, based on the fact that differentiable real valued functions on the real line have 
the Darboux property (\textit{c.f.} Remark~\ref{dense}(i).) Consequently, the main results of this paper apply equally well for the limiting subdifferential.
\end{remark}

%


\subsection{Recovering compact connected sets with nonempty interior}
\label{ss-tapia}

In the current subsection we refine the previous construction to obtain an
everywhere differentiable, compactly supported, $1$-Lipschitz function $f:\mathbb{R}^{d}\rightarrow \mathbb{R}$ such that its subdifferential contains
every closed connected subset of the unit ball with nonempty interior. \smallskip\newline
We shall work with the limiting subdifferential $\partial_{L}f$ which at a given point $x\in \mathbb{R}^N$ consists of all accumulation points of sequences of derivatives 
$\{\nabla f(x_{n})\}_{n\geq 1}$ as $x_{n}\rightarrow x.$ In strong contrast
with the case of strictly differentiable functions (where the only possible limit is $\nabla f(x)$), we show that we can recover all compact connected
sets (even completely irregular fractal-type sets) provided they have nonempty interior.\smallskip \newline
To start, let $C$ be any convex compact set containing $0$ and consider the set
\begin{equation}
\widetilde{\mathcal{K}}_{C}^{0}:=\{K\subset C:~K~\text{is compact connected
and }0\in K\}.  \label{eq:Go 2}
\end{equation}
We first show that, similarly to $\mathcal{K}_{C}^{0},$ the above set can
also be coded on the {Cantor set $\Delta\subset[0,1]$.}
\begin{lemma}
\label{prop_aris} $(\widetilde{\mathcal{K}}_{C}^{0},D_{\mathrm{H}})$ is a {compact metric} space (therefore, it can be seen as continuous surjective image of {the Cantor set $\Delta$).}
\end{lemma}
\noindent \textbf{Proof. }Let us first show that $\widetilde{\mathcal{K}}_{C}^{0}$ is closed in $\mathcal{F}_{C}^{0}$ (see (\ref{eq:F})) for the
Hausdorff distance. To this end, let $\{K_{n}\}_{n}$ be a sequence in $\widetilde{\mathcal{K}}_{C}^{0}$ that converges to a compact set 
$K\in  \mathcal{F}_{C}^{0}.$ It is straightforward to see that $0\in K{\subset C}.$ If $K$ is not connected, then there would exist two nonempty disjoint open subsets 
$U_{1}$ and $U_{2}$ in $\mathbb{R}^{N}$ such that $K^{i}=K\cap U_{i}$ is nonempty, for $i\in \{1,2\}$ and $K=K^{1}\cup K^{2}.$ Then the convergence 
$D_{\mathrm{H}}(K_{n},K)\longrightarrow 0$ forces $K_{n}$ to be disconnected for $n$ sufficiently large, which is a contradiction. {This shows that $(\widetilde{\mathcal{K}}_{C}^{0},D_{\mathrm{H}})$ is a
compact metric space and there exists a continuous surjective function $h$ that maps the Cantor set $\Delta$ onto $\widetilde{\mathcal{K}}_{C}^{0}$ (see \cite[Theorem~4.18]{Kechris})}.\hfill $\Box $

\bigskip 

\noindent {Based on Remark~\ref{rem-gap}(iii),} we can now refine the proof of Lemma~\ref{lem: function Rn} and enhance the conclusion. This is done in the following lemma, whose
proof follows closely the proof of Lemma~\ref{lem: function Rn}. We present a sketch of the proof, highlighting the main changes.\smallskip\newline
Before we proceed, let us recall that a closed set $C\subset \mathbb{R}^{d}$ is called \textit{strictly convex} if for any two distinct points $x,y\in C$, the open segment 
$(x,y)$ joining $x$ and $y$  lies in the interior of $C$. (In particular, a strictly convex set is either singleton or has nonempty interior.)

\begin{lemma} \label{lem: function Rn bis} 
Let $C\subset\mathbb{R}^{N}$ be a convex compact set such that $0\in C$ and $L:=\underset{x\in C}{\max}\{\Vert x \Vert\}.$ Then:\smallskip\newline
(i). There is a differentiable $L$-Lipschitz continuous and compactly supported function $f:\mathbb{R}^{N}\rightarrow\mathbb{R}$ such that:
\begin{equation} \label{eq:salas1 2}
   \text{for every } K\in\widetilde{\mathcal{K}}_{C}^{0}, \text{ there exists } x \in\mathbb{R}^{N} \text{ such that } \partial_L f(x)=K.
\end{equation}
\noindent (ii). Let us further assume that $0\in \mathrm{int}(C)$ and $C$ is strictly convex. Then in addition to the above conclusion we get: 
\begin{equation}\label{eq:salas3 bis}
		\partial_L f(x)\subset C, \text{ for all } x\in\mathbb{R}^{N}
\end{equation}
\end{lemma}
\noindent \textbf{Proof (Sketch).} \textbf{(i)}. Let {$h:\Delta\rightarrow \widetilde{\mathcal{K}}_{C}^{0}$} be a continuous surjective map (which will be used to code the
elements of $\widetilde{\mathcal{K}}_{C}^{0}$). Let $\mathcal{D}=\left\{ d_{n}:\text{ }n\in \mathbb{N}\right\}$ be a countable dense subset of {$\Delta\setminus\{0,1\}$} and consider 
two sequences $\{\alpha _{n}\}_{n}$ and $\{\varepsilon _{n}\}_{n}$
of positive real numbers as in the proof of Lemma~\ref{lem: function Rn}.
In particular, we have $\alpha _{n}>\varepsilon _{n}>0$, for all $n\geq 1,$
$\lim\limits_{n\rightarrow \infty }\alpha _{n}=0$ and $\lim\limits_{n\rightarrow \infty }\varepsilon _{n}/\alpha _{n}=0$.\smallskip\newline
 Define $\{Q_{n}\}_{n}$ by~\eqref{eq:Qin}. Then, the sets $\{B(Q_{n},\varepsilon _{n})\}_{n}$ are
pairwise disjoint and are contained in~$[0,1]^{N}$. Let $\{\gamma _{n}\}_{n}$
be an arbitrary sequence of positive numbers converging to $0$. For every $n\in \mathbb{N}$, since $h(d_{n})$ is totally bounded, there exists a
finite $\gamma _{n}$-net $A_{n}$ of $h(d_{n})$, containing $0$, \textit{i.e.} 
\begin{equation*}
0\in A_{n}\,\subset\, h(d_{n})\,\subset\, \bigcup_{a\in A_{n}}\overline{B}(a,\gamma _{n}).
\end{equation*}
We then define 
\begin{equation}
\widetilde{H}_{n}:=\left( \,A_{n}+\overline{B}(0,2\gamma _{n})\,\right)
\,\bigcap \,\overline{B}(0,L)=\bigcup_{a\in A_{n}}\overline{B}(a,2\gamma_{n})\cap \overline{B}(0,L).  \label{eq:salas3 2}
\end{equation}
Notice that $\widetilde{H}_{n}$ is a finite union of strictly convex sets
and $0\in B\left( 0,\min \left\{ 2\gamma _{n},L\right\} \right) $. 
Moreover, for every $n\in \mathbb{N}$ we have
\begin{equation}
h(d_{n})\,\subset \,\mathrm{int}(\widetilde{H}_{n})\,\subset h(d_{n})+\overline{B}(0,2\gamma _{n}).  \label{eq:set inclusion 2}
\end{equation}
Therefore, $\mathrm{int}(\widetilde{H}_{n})$ is connected, therefore,
according to~\cite[Theorem~8]{BFKL}, for every $n\in \mathbb{N}$, there exists a $\mathcal{C}^{1}$-smooth function $b_{n}:\mathbb{R}^{N}\rightarrow \mathbb{R}$, 
with support in the unit ball, such that $\nabla b_{n}(\mathbb{R}^{N})=\widetilde{H}_{n}$ and $\Vert b_{n}\Vert _{\infty }\leq 1$. 
We set 
\begin{equation}
\phi _{n}(x):=\varepsilon _{n}\cdot b_{n}\Bigl(\frac{x-Q_{n}}{\varepsilon_{n}}\Bigr)  \label{eq:phiin 2}
\end{equation}

\noindent We are ready to define the function $f$ that satisfies our assertion: 
\begin{equation}
\left\{ 
\begin{array}{l}
f:\mathbb{R}^{N}\rightarrow \mathbb{R}\smallskip  \\ 
f(x)=\sum\limits_{n=1}^{\infty }\phi _{n}(x).
\end{array}
\right.   \label{eq:def-sum 2}
\end{equation}
Since $\mathrm{supp}(f)\subset \lbrack 0,1]^{N},$ the function $f$ is
compactly supported.\medskip \newline
\textit{Claim}: For every $K\in \widetilde{\mathcal{K}}_{C}^{0}$ there
exists $x\in [0,1]\times \{0\}^{N-1}$ with $\partial _{L}f(x)=K$.\medskip \newline
\textit{Proof of the Claim}. It follows as in the proof of Lemma~\ref{lem: function Rn} by noticing that~\eqref{eq:set inclusion 2} gives us that 
\begin{equation*}
D_{\mathrm{H}}(\widehat{H}_{n},h(d_{n}))\leq 2\gamma _{n},~\text{for all } n\in \mathbb{N},
\end{equation*}
{and that the Cantor set $\Delta$ is a perfect set.}
\hfill $\lozenge $\newline
Finally, the differentiability of $f$ follows exactly as in the proof of
Lemma~\ref{lem: function Rn}. This completes the proof of (i).\bigskip \newline
\noindent \textbf{(ii).} We now assume that there exists $\lambda >0$ such
that $B(0,\lambda )\subset C$ and that $C$ is strictly convex. To construct
a function $f$ that satisfies (\ref{eq:salas1 2})--(\ref{eq:salas3 bis}), we
replace the definition of $\widetilde{H}_{n}$ in (\ref{eq:salas3 2}) by 
\begin{equation*}
\widetilde{H}_{n}:=\left( A_{n}+\overline{B}(0,2\gamma _{n})\right) \cap
C=\bigcup_{a\in A_{n}}\overline{B}(a,2\gamma _{n})\cap C.
\end{equation*}
Thus, $\widetilde{H}_{n}$ is a finite union of strictly convex sets and $0\in \mathrm{int}(\widetilde{H}_{n})$. Proceeding as before, it easily
follows that $\partial _{L}f(x)\subset C\subset B(0,L),$ for all $x\in \mathbb{R}^{N}$ (in particular $f$ is $L$-Lipschitz) and (\ref{eq:salas1 2})
follows as in (i). \hfill $\Box $ \bigskip \newline
Similarly to the proof of Theorem~\ref{theo: function RN}, we can now use
Lemma~\ref{lem: function Rn bis}  to obtain the existence of a compactly
supported differentiable $1$-Lipschitz function ${f:\mathbb{R}^{N}\rightarrow \mathbb{R}}$ such that the range of its limiting subdifferential $\partial_L f$ contains all compact, connected subsets of the closed unit ball $\overline{B}(0,1)$ with nonempty interior. (Notice that Lemma~\ref{lem: function Rn bis} uses the fact that the Euclidean balls are strictly convex.) Then by a standard argument, already evoked in the beginning of Subsection~\ref{ss-3.2}, see~\eqref{eq:AD} we deduce the following result. \smallskip

\begin{theorem} \label{jortega} There exists a differentiable locally Lipschitz function $f:\mathbb{R}^{N}\rightarrow\mathbb{R}$ such that for every compact,
connected subset $K$ of $\mathbb{R}^N$ with nonempty interior, there exists  $x\in\mathbb{R}^{N}$ such that $\partial_L f(x)=K$.
Moreover, given $\varepsilon>0$, $f$ can be taken to satisfy $\|f\|_{\infty}< \varepsilon$.
\end{theorem}
\noindent Let us mention the following interesting consequence of the above result. Denoting by
\begin{equation*}
\mathrm{gph}(\nabla f):=\{(x,\nabla f(x)):\,x\in \mathbb{R}^{N}\}\,\subset \,
\mathbb{R}^{N}\times \mathbb{R}^{N}
\end{equation*}
the graph of the derivative $\nabla f$ of a differentiable function $f:\mathbb{R}^{N}\rightarrow \mathbb{R}$, we have:
\begin{corollary}
There exists a differentiable locally Lipschitz function $f:\mathbb{R}^{N}\rightarrow\mathbb{R}$ with the property that for every compact, connected subset $K$ of $\mathbb{R}^N$ {with nonempty interior}, there
exists~$\bar x\in\mathbb{R}^{N}$ such that 
\begin{equation}  \label{tapia}
(\bar x, y)\in \overline{\mathrm{gph}(\nabla f)}\Longleftrightarrow ~y\in K.
\end{equation}
\end{corollary}
\noindent This illustrates the gap between mere differentiability versus $\mathcal{C}^1$-smoothness, since in the latter case, only a singleton set $K$ (namely, $K= \{\nabla f(\bar x)\}$) satisfies~\eqref{tapia}.\medskip\newline
\noindent Let us finally notice that Theorem~\ref{jortega} can be seen as a result of \textit{almost exhaustiveness} for the limiting subdifferential of a \textit{differentiable}, locally Lipschitz function. Indeed, Mal\'y~\cite{Maly} established a Darboux-type property for the gradient $\nabla f$ of a differentiable function $f$ in $\mathbb{R}^N$, namely, that  
$$\nabla f(B):=\{\nabla f(x):\,x\in B\}$$ is connected, for any convex body $B$ of $\mathbb{R}^N$. It follows that if $f$ is differentiable and locally Lipschitz, then the above set is bounded and the limiting subdifferential is also given by the formula
\[
\partial_L f(\bar x) = \bigcap_{\varepsilon > 0} \,\mathrm{cl }\left(  \{ \nabla f(x):\, x\in \overline{B}(\bar x, \varepsilon) \} \right).
\]
Therefore $\partial_L f(\bar x)$ contains $\{\nabla f(\bar x)\}$ and is always a compact connected set (as intersection of nested compact connected sets). It follows that the differentiable, locally Lipschitz function $f$ of the statement of Theorem~\ref{jortega} is almost exhaustive for the limiting subdifferential (compare with Definition~\ref{def_exh}). \medskip\newline
\noindent\rule{5cm}{0.5mm} \medskip \\
\noindent\textbf{Acknowledgment.} This work was initiated during a research visit of the second author at TU~Wien (November 2023). This author thanks TU~Wien and the VADOR\ group for hospitality. {The authors are grateful to an anonymous referee for his/her careful reading. This research was funded in whole or in part by the Austrian Science Fund (FWF) [DOI 10.55776/P36344N]. For
open access purposes, the first author has applied a CC BY public copyright license to
any author accepted manuscript version arising from this submission. }

\smallskip


\noindent\rule{5cm}{0.5mm} \bigskip\newline\noindent Aris DANIILIDIS, Sebasti{\'{a}}n
TAPIA-GARC{\'{I}}A

\medskip

\noindent Institute of Statistics and Mathematical Methods in Economics,
E105-04 \newline TU Wien, Wiedner Hauptstra{\ss }e 8, A-1040 Wien\smallskip
\newline\noindent E-mail: \{\texttt{aris.daniilidis,
sebastian.tapia\}@tuwien.ac.at}\newline\noindent
\texttt{https://www.arisdaniilidis.at/}\newline
\texttt{https://sites.google.com/view/sebastian-tapia-garcia}

\medskip

\noindent Research supported by the Austrian Science Fund grant (FWF) \textsc{{DOI 10.55776/}P-36344N}.\newline\vspace{0.2cm}

\noindent\noindent Robert DEVILLE

\smallskip

\noindent Laboratoire Bordelais d'Analyse et Geom\'{e}trie\newline Institut de
Math\'{e}matiques de Bordeaux, Universit\'{e} de Bordeaux 1\newline351 cours
de la Lib\'{e}ration, Talence Cedex 33405, France \medskip

\noindent E-mail: \texttt{Robert.Deville@math.u-bordeaux1.fr}

\end{document}